\documentclass[reqno,twoside,11pt]{article}
\usepackage{amssymb,amsfonts,amsthm,amsmath}
\usepackage{color}
\usepackage{hyperref}%\hypersetup{pagebackref=true}
\usepackage[left=1 in,top=1 in,right=1 in,bottom=1 in]{geometry}
\usepackage{cite}
\newcommand{\beq}{\begin{equation}}
\newcommand{\eeq}{\end{equation}}
\newcommand{\beqs}{\begin{equation*}}
\newcommand{\eeqs}{\end{equation*}}
\newcommand{\ba}{\begin{array}}
\newcommand{\ea}{\end{array}}
\newcommand{\beas}{\begin{eqnarray*}}
\newcommand{\eeas}{\end{eqnarray*}}
\newcommand{\bea}{\begin{eqnarray}}
\newcommand{\eea}{\end{eqnarray}}
\newcommand{\bal}{\begin{align}}
\newcommand{\eal}{\end{align}}

\newcommand{\bals}{\begin{align*}}
\newcommand{\eals}{\end{align*}}

\newcommand{\R}{\ensuremath{\mathbb R}}

\newcommand{\norm}[1]{\| {#1} \|}

\newcommand{\bds}{\begin{displaystyle}}
\newcommand{\eds}{\end{displaystyle}}

\renewcommand{\eqref}[1]{(\ref{#1})}

\def\longequals{\mathbin{=\kern-2pt=}}
\def\eqdef{\stackrel{\rm def}{=}}

\def\varep{\varepsilon}

\def\ddt{\frac{d}{dt}}

\newcommand{\remove}[1]{} %-- ON
\renewcommand{\remove}[1]{#1} % OFF

\newtheorem{theorem}{Theorem}[section]

\newtheorem{lemma}[theorem]{Lemma}
\newtheorem{corollary}[theorem]{Corollary}
\newtheorem{proposition}[theorem]{Proposition}

\newtheorem{remark}[theorem]{\bf{Remark}}

\theoremstyle{remark}

\def\midx{\mu}
\def\ellz{{\ell_Z}}
\def\ellf{{\ell_f}}
\def\ellb{{\ell_B}}
\def\myclearpage{}

\definecolor{darkred}{rgb}{.70,.12,.20}

\definecolor{darkgreen}{rgb}{.20,.52,.14}

\numberwithin{equation}{section}

\title{Doubly nonlinear parabolic equations for a general class of Forchheimer gas flows in porous media}
\author{Emine Celik$^\dag$, Luan Hoang$^\dag$, and Thinh Kieu$^\ddag$}

\begin{document}
\date{\today}
\maketitle 
\begin{center}
\textit{$^\dag$Department of Mathematics and Statistics, Texas Tech University, Box 41042, Lubbock, TX 79409--1042, U. S. A.} \\
\textit{$^\ddag$Department of Mathematics, University of North Georgia, Gainesville Campus, 3820 Mundy Mill Rd., Oakwood, GA 30566, U. S. A.}\\
Email addresses: \texttt{emine.celik@ttu.edu, luan.hoang@ttu.edu, thinh.kieu@ung.edu}\\
%\bigskip{\large \today}
\end{center}
%\setcounter{equation}{0}
%-----------------------------------------------------------------------------------------------------------------
\begin{abstract} 
This paper is focused on the generalized Forchheimer flows of compressible fluids in porous media. The gravity effect and other general nonlinear forms of   the source terms and boundary fluxes are integrated into the model. 
It covers isentropic gas flows, ideal gases and slightly compressible fluids. We derive a doubly nonlinear parabolic equation for the so-called pseudo-pressure, and study the corresponding initial boundary value problem. The maximum estimates of the solution are established by using suitable trace theorem and adapting appropriately the Moser's iteration. The gradient estimates are obtained under a theoretical condition which, indeed, is relevant to the fluid flows in applications. 
\end{abstract}
%-----------------------------------------------------------------------------------------------------------------

%\tableofcontents 
%   \newpage \pagenumbering{arabic}
   %\tableofcontents 

\pagestyle{myheadings}\markboth{E. Celik, L. Hoang, and T. Kieu}
{A General Class of Forchheimer Gas Flows in Porous Media}

%%%%%%%%%%%%%%%%%%%%%%%%%%%%%%%%%%%%%%%%%%%
\myclearpage    
\section{Introduction}\label{Intro}
%%%%%%%%%%%%%%%%%%%%%%%%%%%%%%%%%%%%%%%%%%%
We consider fluid flows in porous media with pressure $p$, density $\rho$, velocity $v$, and absolute viscosity $\mu$.
The media has permeability $k>0$ and porosity $\phi\in(0,1)$. For hydrodynamics in  porous media, the following Darcy's equation  is usually used as a default law
\beq\label{Darcy}
-\nabla p = \frac {\mu}{k} v.
\eeq

However, even  Darcy himself \cite{Darcybook} noted that there were deviations from the linear equation \eqref{Darcy}. For instance,  in the case the Reynolds number is large or fluids in fractured media, \eqref{Darcy}  becomes inaccurate in describing the fluid dynamics.  Many work have been devoted to developing alternative nonlinear models to  Darcy's law, see e.g. \cite{BearBook}.
%Nonlinear equations, therefore, are needed,  see \cite{MuskatBook,BearBook}. 
%Even in early works,  Darcy \cite{Darcybook} and Dupuit \cite{Dupuit1857} already acknowledged the deviations from equation \eqref{Darcy}.
Forchheimer, in \cite{Forchh1901,ForchheimerBook}, established the following three models:  
\beq\label{2term}
-\nabla p=av+b|v|v,
\eeq
\beq\label{3term}
-\nabla p=av+b|v|v+c |v|^2 v,
\eeq
\beq\label{power}
-\nabla p=av+d|v|^{m-1}v,\quad\text{for some real number } m\in(1, 2).
\eeq
The numbers $a,b,c,d$ above are empirical positive constants.
The equations \eqref{2term}, \eqref{3term} and \eqref{power} are usually referred to as Forchheimer's two-term, three-term and power laws, respectively.
For more models and discussions, see \cite{MuskatBook,Ward64,BearBook,NieldBook,StraughanBook} and references therein. 

From the mathematical point of view, the Darcy flows,  under the umbrella ``porous medium equations", have been analyzed intensively  since 1960s,  see e.g. \cite{VazquezPorousBook} and a large number of references cited there. 
In contrast, the mathematics of  Forchheimer flows and their variations has attracted much less attention.
Moreover, the existing papers on this topic mainly treat incompressible fluids, leaving compressible ones barely studied,
see \cite{StraughanBook} and references therein. The current paper aims to explore Forchheimer flows  of compressible fluids using analytical techniques from partial differential equation (PDE) theory.

The Forchheimer  equations \eqref{2term}, \eqref{3term}, \eqref{power} are extended to more general form 
\beq\label{gF}
-\nabla p =\sum_{i=0}^N a_i |v|^{\alpha_i}v,
\eeq  
where $a_i$'s are positive constants.  
Equation \eqref{gF} is called the generalized Forchheimer equation.
It is used to unify the models \eqref{2term}, \eqref{3term}, \eqref{power}, and as a framework for interpretation of different experimental or field data.
It is analyzed numerically in \cite{Doug1993,Park2005, Kieu1},
theoretically in \cite{ABHI1,HI2,HIKS1,HKP1,HK1,HK2} for single-phase flows, and also in \cite{HIK1,HIK2} for two-phase flows.

For compressible fluids, especially gases, the dependence of coefficients $a_i$'s on the density $\rho$ is essential.
By using dimension analysis, Muskat \cite{MuskatBook} and then Ward \cite{Ward64} proposed the following equation for both laminar and turbulent flows in porous media:
\beq\label{W}
-\nabla p =f(v^\alpha k^{\frac {\alpha-3} 2} \rho^{\alpha-1} \mu^{2-\alpha}),\text{ where  $f$ is a function of one variable.}
\eeq 

Using this, Ward \cite{Ward64} established from experimental data that
\beq\label{FW} 
-\nabla p=\frac{\mu}{k} v+c_F\frac{\rho}{\sqrt k}|v|v\quad \text{with }c_F>0.
\eeq
This model is widely accepted as the standard form of the Forchheimer's two term-law.

Based on the arguments by  Muskat and Ward, we proposed in \cite{CHK1} the following adaptation for \eqref{gF}
\beq\label{FM}
-\nabla p= \sum_{i=0}^N a_i \rho^{\alpha_i} |v|^{\alpha_i} v,
 \eeq
where $N\ge 1$, $\alpha_0=0<\alpha_1<\ldots<\alpha_N$ are real numbers, the coefficients $a_0, \ldots, a_N$ are positive.  
This equation covers the two-term case \eqref{FW}, and the focus is  the dependence on the density, but not viscosity and permeability.

Equation \eqref{FM}, however, does not take into account the gravity. Because of the nonlinear density-dependence of the model, any addition of new density terms may complicate the analysis as we will see below.
Nonetheless, the gravity can be integrated into \eqref{FM} by replacing  $(-\nabla p)$ with $-\nabla p + \rho \vec g$, where $\vec g$ is the constant gravitational field. Therefore, we  consider 
\beq\label{anisot}
 \sum_{i=0}^N a_i \rho^{\alpha_i} |v|^{\alpha_i} v  =-\nabla p + \rho \vec g.
\eeq
%with positive coefficients $a_i$'s.

%***Here, the viscosity and permeability are considered constant and we do not specify the dependence of $a_i$'s on them.
%Our mathematical exposition below will allow all $\alpha_i\ge 0$  in \eqref{FM}. 
%In practice, we can simply take $\alpha_N\le 2$ in \eqref{FM} or use the popular model \eqref{FW}.
%Even in these cases,  the results obtained in this paper are new.***

Denote by $g:\mathbb{R}^+\rightarrow\mathbb{R}^+$ a generalized polynomial with positive coefficients defined by
\beq\label{eq2}
g(s)=a_0s^{\alpha_0} + a_1s^{\alpha_1}+\cdots +a_Ns^{\alpha_N}\quad\text{for } s\ge 0,
\eeq 
with $a_0,a_1,\ldots,a_N>0$.
Then \eqref{anisot} can be rewritten as
\beq\label{eq0}
g(\rho|v|) v  =-\nabla p + \rho \vec g,
\eeq

With this rule,  we now follow the method in \cite{CHK1} to derive the basic PDE and its corresponding boundary condition.
Multiplying both sides of \eqref{eq0} by $\rho$, we obtain 
 \beq\label{eq1}
 g(|\rho v|) \rho v   =- \rho\nabla p +\rho^2 \vec g.
 \eeq

 Solving for $\rho v$ from \eqref{eq1} gives
 \beq\label{solve} 
\rho v=-K(|\rho\nabla p -\rho^2 \vec g|)(\rho\nabla p -\rho^2 \vec g), 
\eeq
where the function $K: \mathbb{R}^+\rightarrow\mathbb{R}^+$ is defined for $\xi\ge 0$ by
\beq \label{Kdef}
K(\xi)=\frac{1}{g(s(\xi))},
\eeq
 with  $s(\xi)=s$  being the unique non-negative solution of  $s g(s)=\xi.$

Relation \eqref{solve} will be combined with other equations of fluid mechanics. The first is the continuity equation
\beq\label{con-law}
\phi\rho_t+{\rm div }(\rho v)=F,\eeq
where  the porosity$\phi$  is a constant in $(0,1)$, the \textit{source term} $F$ counts for the rate of net mass production or loss due to any source and/or sink in the media.

From \eqref{solve} and \eqref{con-law} follows
\beq\label{genr}
\phi\rho_t={\rm div }(K(|\rho\nabla p -\rho^2 \vec g|)(\rho\nabla p -\rho^2 \vec g))+F.
\eeq

We consider below scenarios of isentropic gas flows, ideal gases  and slightly compressible fluids.

\textbf{Isentropic gas flows.} In this case
\beq\label{isen}
p=\bar c\rho^\gamma\quad \text{for some constants }  \bar c,\gamma>0.
\eeq
Here, $\gamma$ is the specific heat ratio. 
Note that $\rho\nabla p=\nabla(\bar c\gamma\rho^{\gamma+1}/(\gamma+1))$.
Hence by letting
 \beq \label{urho}
 u=\frac{\bar c\gamma\rho^{\gamma+1}}{\gamma+1}=\frac{\gamma p^\frac{\gamma+1}\gamma}{\bar c^\frac1\gamma(\gamma+1)},
 \eeq
we rewrite \eqref{genr} as
\beq\label{maineq2}
\phi c^{1/2} (u^\lambda)_t= \nabla\cdot(K(|\nabla u-c u^\ell \vec g|)(\nabla u-c u^\ell \vec g)) +F,
\eeq
where
 \beq\label{gaspara}
 \lambda=\frac1{\gamma+1}\in (0,1),\quad \ell=2\lambda \quad \text{and}\quad c=\left(\frac{\gamma+1}{\bar c\gamma}\right)^\ell.
\eeq
The new quantity $u$ in \eqref{urho} is essentially a pseudo-pressure. 
% \beq
%\rho=(\frac{\gamma+1}{c\gamma})^\lambda u^\lambda=c_2^{1/2}u^\lambda, \quad \lambda=\frac1{\gamma+1}\in (0,1].
%\eeq

\textbf{Ideal gases.} The  equation of state is
\beq\label{ideal}
p=\bar c\rho\quad \text{for some constant }  \bar c>0.
\eeq
We can consider \eqref{ideal} as a special case of \eqref{isen} with $\gamma=1$, and derive the same equation \eqref{maineq2} with $\lambda=1/2$.
In this case, the pseudo-pressure $u\sim p^2$ which is used commonly in engineering problems.

Although we mainly focus on gases, we also present here, in a unified way,  the slightly compressible fluids which is important in petroleum engineering.

\textbf{Slightly compressible fluids.} The equation of state is
\beqs
\frac 1\rho \frac{d\rho}{dp}=\frac 1\kappa=const.>0.
\eeqs
Then $ \rho \nabla p=\kappa\nabla \rho$, and by letting 
\beq\label{ukap} u=\kappa \rho,\quad  \ell=2 \quad\text{and}\quad  c=1/\kappa^2,
\eeq
we obtain the same equation \eqref{maineq2} with $\lambda=1$.

\medskip
For all three cases, by scaling the time variable, we can always assume the multiplying factor on the left-hand side of \eqref{maineq2}  to be $1$. In summary, we have derived
\beq\label{maineq}
(u^\lambda)_t= \nabla\cdot(K(|\nabla u-c u^\ell \vec g|)(\nabla u-c u^\ell \vec g)) +F
\eeq
with constants  $\lambda\in(0,1]$, $\ell=2\lambda$, $c>0$, and function $F(x,t)$ being rescaled appropriately.

\textbf{Boundary condition.} 
We will study the problem in a bounded domain $U$ with outward normal vector $\vec \nu$ on the boundary.
We  consider the volumetric flux condition 
$$v\cdot \vec\nu=\psi\text{ on }\partial U.$$
This gives  $\rho v\cdot \vec\nu =\psi \rho$, hence, together with \eqref{solve} and \eqref{urho} or \eqref{ukap}  yields 
\beq\label{bc0}
-K(|\nabla u+c u^\ell \vec g|)(\nabla u+c u^\ell \vec g)\cdot \vec\nu=c^{1/2}\psi u^\lambda .
\eeq

\textbf{General formulation and the  initial boundary value problem (IBVP).} 
%Our equations of interest are \eqref{maineq2} and \eqref{bc0}.
%We rewrite and formulate them into the following general IBVP
Although problem \eqref{maineq} and \eqref{bc0} is our motivation, in this mathematical investigation, we  consider a more general class of equations and boundary conditions, namely,
\beq\label{upb}
\begin{cases}
\begin{displaystyle}\frac{\partial (u^{\lambda})}{\partial t} \end{displaystyle}
= \nabla \cdot (K (|\nabla u+Z(u)|)(\nabla u+Z(u))  )+f(x,t,u) &\text{ on }  U\times (0,\infty),\\
u(x,0)=u_0(x) &\text{ on } U,\\
K (|\nabla u+Z(u)|)(\nabla u+Z(u))\cdot \vec\nu=B(x,t,u)  &\text{ on } \Gamma \times(0,\infty),
\end{cases}
\eeq 
 where  $Z(u)$ is a function from $[0,\infty)$ to $\R^n$, $B(x,t,u)$ is a function from $\Gamma\times[0,\infty)\times[0,\infty)$ to $\R$ and $f(x,t,u)$ is a function from $U\times[0,\infty)\times[0,\infty)$ to $\R$. 
  
In \eqref{upb} the source term $f(x,t,u)$ now can depend on $u$, and the boundary term $B(x,t,u)$ can be more general than $\psi(x,t)u^\lambda$. For our analysis, they still need some growth conditions.
 
 \bigskip
\noindent\textbf{Assumption (A1).} Throughout this paper, we assume that functions $Z(u):[0,\infty)\to \R^n$, $B(x,t,u):\Gamma\times[0,\infty)\times[0,\infty)\to\R$ and $f(x,t,u):U\times[0,\infty)\times[0,\infty)\to \R$ satisfy
\beq\label{Zprop}
|Z(u)|\le d_0 u^\ellz,
%d_0^{-1} u^\ellz\le |Z(u)|\le d_0 u^\ellz,
\eeq
\beq \label{Bprop}
 B(x,t,u)\le \varphi_1(x,t)+\varphi_2(x,t)u^\ellb,
\eeq
\beq\label{fprop}
f(x,t,u)\le f_1(x,t)+f_2(x,t)u^\ellf
\eeq
with constants $d_0,\ellz>0$, $\ellf,\ellb\ge 0$, and functions $\varphi_1,\varphi_2,f_1,f_2\ge 0$.

\bigskip
Note that the growths with respect to $u$ in Assumption (A1) are arbitrary. This is different from other existing papers when the exponents $\ellz$, $\ellf$, $\ellb$ are restricted to suit certain Sobolev embedding or trace theorems.

The problem \eqref{upb} is significantly more general than the one in our previous work \cite{CHK1}. It can arise from  other complex, nonlinear models of fluid flows in porous media.
Indeed, a similar PDE with $p$-Laplacian structure was derived in \cite{DiazThelin1994}  for water vapor.
Here, we showed that it comes naturally from the Forchheimer equations and is formulated for the pseudo-pressure $u$ instead. Due to more complicated flows \eqref{eq0},  the resulting function $K(\xi)$ is non-homogeneous compared to the homogeneous one in \cite{DiazThelin1994}.

Regarding the PDE in \eqref{upb}, it is doubly nonlinear in both $u$ and $\nabla u$.
(For the theory of equations of this type, see the monograph   \cite{IvanovBook82}, review paper \cite{IvanovReg1997} and e.g.  \cite{kinnunen2007local,Tsutsumi1988,ManfrediVespri,Vespri1992,Alt1983}.)
Moreover, our equation contains lower order terms of arbitrary growths in $u$. 
Therefore, it is not clear whether $L^\infty$-estimates are possible.
In addition, the boundary condition is time-dependent, of non-linear Robin type and also has arbitrary growth rate in $u$.
Thus, the boundary contribution is not trivial and, thanks to the nonlinearity, we cannot shift the solution by subtracting the boundary data. 
For vanishing Dirichlet boundary condition, the work \cite{DiazThelin1994}  uses maximum principle which is not applicable to our problem.
%, see also  Ivanov, Jager \cite{ivanov2000existence} for gradient estimates. 
Also, both  \cite{DiazThelin1994,Tsutsumi1988} imposes the $L^\infty$-requirement for the initial data.
%, with much simpler boundary conditions, i.e., vanishing boundary values,
In contrast, we use Moser's iteration \cite{moser1971pointwise},  hence can deal with more complex equation and boundary condition. 
Furthermore, we only require  the initial data belonging to a certain $L^\alpha$-space with finite $\alpha>0$, and derive the $L^\infty$-estimates of the solution for positive time. We note that, although ours are a priori estimates, they are crucial in establishing, via regularization and approximation,  the existence results (see \cite{IvanovBook82}).

The paper is organized as follows.
In section \ref{Prelim}, we recall needed trace theorem and Poincar\'e-Sobolev inequalities. In particular, the inequalities in Corollary \ref{ecor} are formulated to suit the nonlinear diffusion, and later treatment of the general  nonlinearity in the source term and Robin boundary condition. 
%--------------------------------------------
In section \ref{lasec},  we establish the $L^\alpha$-estimate of a solution $u(x,t)$ of \eqref{upb} for any finite  $\alpha>0$  in terms of the initial and boundary data. 
%--------------------------------------------
In section \ref{maxsec}, we derive an estimate for spatially global $L^\infty$-norm of $u(x,t)$ in Theorem \ref{thm45} by adapting Moser's iteration. The sequence of exponents in the iteration are constructed based on the nonlinearity of the boundary condition and the source term. We note that the global $L^\infty$-norm of $u$ is bounded by \eqref{Li1} in Theorem \ref{LinfU} which is ``quasi-homogeneous'' in its $L^\beta_{x,t}$-norm for some $\beta>0$. This extends previous results in \cite{Surnachev2012,CHK1}. 
%--------------------------------------------
In section \ref{gradsec}, we establish  $L^{2-a}$-estimates for the gradient of $u(x,t)$ in Theorem \ref{thm52}. Even for this simple norm, it is non-trivial due to the arbitrary growth in the nonlinear Robin boundary condition. It is obtained under condition \eqref{lamel}. This mathematical requirement turns out to be naturally satisfied for the original problem \eqref{maineq} and known gases such as those in the data book \cite{GasBook2001}. It modestly shows the relevance of our mathematical analysis.

%%%%%%%%%%%%%%%%%%%%%%%%%%%%%%%%%%%%%%%%%%%%%%%%%%%%%%%%%%%%%%%%%%%%%
\myclearpage    
\section{Basic inequalities}\label{Prelim}
%%%%%%%%%%%%%%%%%%%%%%%%%%%%%%%%%%%%%%%%%%%%%%%%%%%%%%%%%%%%%%%%%%%%%

First, we recall elementary inequalities that will be used frequently. 
Let $x,y\ge 0$. Then
%\beq\label{ee1}
%(x+y)^p\le 2^p(x^p+y^p)\quad  \text{for all }  p>0,
%\eeq
\beq\label{ee2}
(x+y)^p\le x^p+y^p\quad  \text{for all } 0<p\le 1,
\eeq
\beq\label{ee3}
(x+y)^p\le 2^{p-1}(x^p+y^p)\quad  \text{for all }  p\ge 1,
\eeq
\beq\label{ee4}
x^\beta \le x^\alpha+x^\gamma\quad \text{for all } 0\le \alpha\le \beta\le\gamma,
\eeq
particularly,
\beq\label{ee5}
x^\beta \le 1+x^\gamma \quad \text{for all } 0\le \beta\le\gamma.
\eeq

Next, we recall particular  Poinc\'are-Sobolev inequality and trace theorem. 
 
The classical trace theorem: If a function $u(x)$ belongs to $W^{1,1}(U)$ then
\beq\label{W11}
\int_\Gamma|u|d\sigma\le c_1\int_U |u|dx +c_2\int_U |\nabla u|dx,
\eeq
where constants $c_1,c_2>0$ depend on $U$.
 
\begin{lemma}[\cite{CHK1}, Lemma 2.1]\label{egentrace}
In the following statements, $u(x)$ is a function  defined on $U$.
\begin{enumerate}
%%%%
\item[\rm(i)] If $\alpha \ge s\ge 0$, $\alpha\ge 1$, and $p>1$, then for any $|u|^\alpha\in W^{1,1}(U)$ and $\varepsilon>0$ one has
\beq\label{etrace10}
\int_\Gamma |u|^\alpha d\sigma 
\le \varepsilon \int_U |u|^{\alpha-s}|\nabla u|^p dx + c_1 \int_U |u|^\alpha dx  + (c_2 \alpha)^\frac p{p-1} \varepsilon^{-\frac 1{p-1}} \int_U |u|^{\alpha+\frac{s-p}{p-1}} dx.
\eeq
%where   $c_1,c_2>0$ are constants depending on $U$, but not on $u(x),\alpha,s,p$.
%%%%
\item[\rm(ii)] If  $n>p>1$, $r>0$, $\alpha\ge s\ge 0$,  $\alpha\ge \frac{p-s}{p-1}$, and $\alpha>\frac{n(r+s-p)}{p}$, then for any $\varepsilon>0$ one has
\begin{multline}\label{eS}
\int_U |u|^{\alpha+r}dx
\le \varepsilon \int_U |u|^{\alpha-s}|\nabla u|^p dx+ \varepsilon^{-\frac \theta{1-\theta}} 2^\frac{\theta(\alpha-s+p)}{1-\theta}(\bar c_3 m)^{\frac{\theta p}{1-\theta}} \|u\|_{L^\alpha}^{\alpha+\midx}\\  
 +2^{\theta (\alpha-s+p)}\bar c_4^{\theta p}|U|^\frac{\theta(\alpha(p-1)+s-p)}{\alpha}\|u\|_{L^\alpha}^{\alpha+r},
\end{multline}
for all $|u|^m\in W^{1,p}(U)$, where  
\beq\label{emdef}
 m=\frac{\alpha-s+p}p,\quad
\theta=\frac{r}{(\alpha p/n)+ p-s},\quad
\midx=\frac{r+\theta(s-p)}{1-\theta},
\eeq
and positive constants $\bar c_3=\bar c_3(U,p)$ and $\bar c_4=\bar c_4(U,p)$ depend on $U,p$, but not on $u(x)$, $\alpha$, $s$.
%%%
\end{enumerate}
\end{lemma}
%==================================================================%

The following remarks on Lemma \ref{egentrace} are in order.

(a) We can calculate $r$ in terms of $\theta$ by the second formula in \eqref{emdef}, and rewrite the power $\alpha+\midx$ in \eqref{eS} as
\beq\label{nmu}
\alpha+\midx=\alpha+\Big[\theta(\frac{\alpha p}{n}+(p-s)) +\theta(s-p)\Big]\frac1{1-\theta}
=\alpha\Big(1+\frac pn\cdot\frac{\theta}{1-\theta}\Big).
\eeq

Assume $p\le s$, then $m\le \alpha$ and we can rewrite \eqref{eS} as 
\beq\label{eS10}
\int_U |u|^{\alpha+r}dx
\le \varepsilon \int_U |u|^{\alpha-s}|\nabla u|^p dx+\bar D_1\|u\|_{L^\alpha}^{\alpha+r}+ \varepsilon^{-\frac \theta{1-\theta}} \bar D_2 \|u\|_{L^\alpha}^{\alpha(1+\frac pn\cdot\frac{\theta}{1-\theta})},
\eeq
where 
\beqs
\bar D_1= [2^{\alpha-s+p}\bar c_4^p|U|^\frac{\alpha(p-1)+s-p}{\alpha}]^{\theta},\quad
\bar D_2=  [2^{\alpha-s+p}(\bar c_3 \alpha )^p]^\frac{\theta}{1-\theta}.
\eeqs

(b) For any $r\ge 0$, $\alpha\ge s\ge 0$, $\alpha+r\ge 1$, applying  \eqref{etrace10} for $\alpha\to\alpha+r$ and $s\to s+r$ yields
\beq\label{tG}
\int_\Gamma |u|^{\alpha+r} d\sigma 
\le \varepsilon \int_U |u|^{\alpha-s}|\nabla u|^p dx + c_1 \int_U |u|^{\alpha+r} dx  + (c_2 (\alpha+r))^\frac p{p-1} \varepsilon^{-\frac 1{p-1}} \int_U |u|^{\alpha+\hat r} dx,
\eeq
where
\beqs
\hat r= r+\frac{r+s-p}{p-1}=\frac{pr+s-p}{p-1}.
\eeqs

%======================================================%

\medskip
In our particular case, we have following corollary.

\begin{corollary}\label{ecor}
Assume $1>a>\delta\ge 0$, $\alpha\ge 2-\delta$, and $|u|^\alpha\in W^{1,1}(U)$. Let $r>0$. 
 \begin{enumerate}
\item[\rm(i)] For any $\varepsilon>0$ one has
\begin{multline}\label{t7}
\int_\Gamma |u|^{\alpha+r} d\sigma\le \varep\int_U |u|^{\alpha-2+\delta}|\nabla u|^{2-a} dx + c_1 \int_U |u|^{\alpha+r} dx \\
 + (c_2 (\alpha+r))^\frac {2-a}{1-a} \varepsilon^{-\frac 1{1-a}} \int_U |u|^{\alpha+\frac{(2-a)r+a-\delta}{1-a}} dx.
\end{multline}

\item[\rm(ii)] If $\alpha>\frac{n(r+a-\delta)}{2-a}$, then for any $\varepsilon>0$ one has
\beq\label{t6}
\int_U |u|^{\alpha+r}dx
\le \varepsilon \int_U |u|^{\alpha-2+\delta}|\nabla u|^{2-a} dx+ D_1\|u\|_{L^\alpha}^{\alpha+r}+D_2\varepsilon^{-\frac{\theta}{1-\theta}}\|u\|_{L^\alpha}^{\alpha\big(1+\frac{2-a}{n}\cdot\frac{\theta}{1-\theta}\big)},
\eeq
where 
\beqs
\begin{aligned}
\theta&=\frac{r}{\alpha(2-a)/n+\delta-a},\\
D_1&=(2^{\alpha+\delta-a}c_4^{2-a}|U|^\frac{\alpha(1-a)+a-\delta}{\alpha})^{\theta} ,
\quad D_2=(2^{\alpha+\delta-a}(c_3\alpha)^{2-a})^\frac{\theta}{1-\theta}
\end{aligned}
\eeqs
with $c_3=\bar c_3(U,2-a)>0$ and $c_4=\bar c_4(U,2-a)>0$.
\end{enumerate}
\end{corollary}
%=====================
\begin{proof}
Let $p=2-a$ and $s=2-\delta$, noticing that $p< s$.
Then inequality \eqref{t7} follows \eqref{tG}, and inequality \eqref{t6} follows \eqref{eS10}.
%\alert{$c_3$ and $c_4$ are fixed with $2-a$.}
\end{proof}

We also recall a particular multiplicative parabolic Sobolev inequality.

%=================================
\begin{lemma}[\cite{CHK1}, Lemma 2.3]\label{Parabwk}
Assume $1>a>\delta\ge 0$,
\beq \label{alpcond}
\alpha\ge 2-\delta\quad{and}\quad \alpha > \alpha_*\eqdef \frac{n(a-\delta)}{2-a}.
\eeq %(comes from $\kappa>1$)
If $T>0$, then
\begin{align}
\Big(\int_0^T \int_U |u|^{\kappa\alpha}dx dt \Big)^{\frac 1{\kappa\alpha}}
&\le (c_5\alpha^{2-a})^\frac1{\kappa\alpha} \Big(  \int_0^T\int_U |u|^{\alpha+\delta-a}dxdt+\int_0^T\int_U |u|^{\alpha+\delta-2}|\nabla u|^{2-a}dxdt \Big)^{\frac {\tilde \theta}{\alpha+\delta-a}}\notag \\
&\quad \cdot\sup_{t\in[0,T]}\Big( \int_U|u(x,t)|^{\alpha}dx \Big)^{\frac{1-\tilde \theta}{\alpha}},\label{parabwki}
\end{align}
where  $c_5\ge 1$ is independent of $\alpha$ and T, and 
\beq \label{kappadef}
\tilde \theta=\tilde \theta_\alpha \eqdef \frac 1{1+\frac{\alpha(2-a)}{n(\alpha+\delta-a)}},\quad 
\kappa=\kappa(\alpha)  \eqdef 1+\frac{2-a}n-\frac{a-\delta}{\alpha}=1+(a-\delta)\Big(\frac1{\alpha_*}-\frac1{\alpha}\Big).
\eeq
\end{lemma}
%=================================

%%%%%%%%%%%%%%%%%%%%%%%%%%%%%%%%%%%%%%%%%%%%%%%%%%%%%%%%%
\myclearpage    
\section{Estimates of the Lebesgue norms }\label{lasec}
%%%%%%%%%%%%%%%%%%%%%%%%%%%%%%%%%%%%%%%%%%%%%%%%%%%%%%%%%

%The class of functions $g(s)$ as in \eqref{FM} is denoted by $FP(N,\vec{\alpha})$ which is the abbreviation of ``Forchheimer polynomials". When the function $g$ in \eqref{eq1} is one of the $g(s)$ in \eqref{eq2}, it is referred to as the Forchheimer polynomial.   

%Let $g=g(s,\vec{a})$ in $FP(N,\vec{\alpha})$. 

From now on, we fix a function $g(s)$ in \eqref{eq2} and \eqref{eq0}. 
Denote  
\beq\label{eq9}
a=\frac{\alpha_N}{\alpha_N+1}\in (0,1).
\eeq
This number will be used in our calculations throughout.

The function $K(\xi)$ defined by \eqref{Kdef} has the following properties: it is decreasing in $\xi$,  maps $[0,\infty)$ onto $(0,\frac{1}{a_0}]$    and
  \begin{align}
\frac{d_1}{(1+\xi)^a}\le K(\xi)\le \frac{d_2}{(1+\xi)^a},\label{i:ineq1} \\ 
  d_3(\xi^{2-a}-1)\le K(\xi)\xi^2\le d_2\xi^{2-a}, \label{i:ineq2}
  \end{align}
  where $d_1, d_2, d_3$ are positive constants depending on $\alpha_i$'s and $ a_i$'s of the function $g(s)$, see \cite{ABHI1}.
Let $n=2,3,\ldots,$ and $U$ be the a bounded, open, connected subset of $\mathbb{R}^n$ with $C^2$ boundary $\Gamma=\partial U$. In applications, $n=2$ or $3$, but we treat all $n\ge 2$ in this paper.

Hereafter, $u(x,t)\ge 0$ is a solution of IBVP \eqref{upb}.

Define $\delta=1-\lambda\in [0,1)$. Throughout the paper, we assume 
\beq\label{ad}
a>\delta.
\eeq

This assumption is to avoid too many possible cases for our estimates. The case $a\le \delta$ can be treated similarly and, in fact, is easier to deal with.

Then  $\alpha_*= n(a-\delta)/(2-a)$ in \eqref{alpcond} is  a fixed positive number.
In order to describe our results, we introduce some constants and quantities.
Define
\beq\label{k*def}
\kappa_B= 1+\frac{1-a}n \quad\text{and}\quad \kappa_f= 1+\frac{2-a}n.
\eeq 

Let $p_1$, $p_2$, $p_3$, $p_4$ be fixed numbers such that
\beq\label{pscond}
\begin{split}
1<p_1,p_2<\kappa_B\quad\text{and}\quad
1<p_3,p_4<\kappa_f.
\end{split}
\eeq

For $i=1,2,3,4$ let $q_i$  be the conjugate exponent of $p_i$, that is, $1/p_i+1/q_i=1$. Then \eqref{pscond} is equivalent to
\beq\label{qscond}
1+\frac n{1-a}<q_1,q_2<\infty\quad\text{and}\quad
1+\frac n {2-a}<q_3,q_4<\infty.
\eeq

A key quantity in expressing our estimates is 
\beq\label{Upsilondef}
\Upsilon(t)= \|\varphi_1(t)\|_{L^{q_1}(\Gamma)}^{q_1}+\|\varphi_2(t)\|_{L^{q_2}(\Gamma)}^{q_2} + \|f_1(t)\|_{L^{q_3}(U)}^{q_3}+ \|f_2(t)\|_{L^{q_4}(U)}^{q_4} \quad \text{for }t\ge 0.
\eeq

%$$\varrho,\varpi,\vartheta,\zeta,\varsigma$$

We also denote by $\eta_0$ a positive number defined by  
\beq\label{eta1}
\begin{aligned}
\eta_0
=\max &\Big\{  q_1\lambda, q_2(\lambda-\ellb),  n(\ellz-1), 
\frac{-p_1\lambda+a-\delta}{\kappa_f-p_1},
\frac{p_2(-\lambda+\ellb)+a-\delta}{\kappa_f-p_2},\\
&\quad \frac{-p_1\lambda+(a-\delta)}{\kappa_B-p_1},
\frac{p_2(-\lambda+\ellb)
+a-\delta}{\kappa_B-p_2},
\frac{-p_3\lambda+a-\delta}{\kappa_f-p_3},
\frac{p_4(-\lambda+\ellf)+a-\delta}{\kappa_f-p_4}
\Big\}.
\end{aligned}
\eeq

In the following, we  focus on estimating the $L^\alpha$-norm on $U$ of the solution  $u(x,t)$, for any given $\alpha>0$ and $t>0$. We start with a differential inequality for $\|u(t)\|_{L^\alpha(U)}^\alpha$.

%==========================================================%
\begin{lemma}\label{lalem}
 Assume 
 \beq \label{al3}
 \alpha\ge \max\Big\{2,\frac{na}{1-a}\Big\}\quad\text{and}\quad \alpha>\eta_0.
 \eeq
 
For $t>0$, one  has
\beq\label{est10}
%\begin{aligned}
\frac{d }{dt}\int_U u(x,t)^{\alpha} dx+ \int_U |\nabla u(x,t)|^{2-a}u(x,t)^{\alpha-\lambda-1}dx\le C_{0}\cdot \Big(\|u(t)\|_{L^\alpha(U)}^{\nu_1} +\|u(t)\|_{L^\alpha(U)}^{\nu_2} +\Upsilon(t)\Big),
%\end{aligned}
\eeq
where $C_{0}=C_{0}(\alpha)>0$, 
\beqs
\nu_1=\nu_1(\alpha)\eqdef \alpha-\lambda-1\quad\text{and}\quad \nu_2=\nu_2(\alpha)\eqdef \alpha\Big(1+\frac{2-a}{n}\cdot \frac\theta{1-\theta}\Big)
\eeqs
with 
\beqs
\theta=\theta_\alpha\eqdef \frac{\nu_3-\alpha}{\alpha(2-a)/n+\delta-a}\quad \text{for $\nu_3=\nu_3(\alpha)$ defined by \eqref{mu3} below}.
\eeqs
%and
%\beq\label{Upsilondef}
%\Upsilon(t)= \|\varphi_1(t)\|_{L^{q_1}(\Gamma)}^{q_1}+\|\varphi_2(t)\|_{L^{q_2}(\Gamma)}^{q_2} + \|f_1(t)\|_{L^{q_3}(U)}^{q_3}+ \|f_2(t)\|_{L^{q_4}(U)}^{q_4},
%\eeq
\end{lemma}
%==========================================================%
\begin{proof}
Multiplying the first equation in \eqref{upb} by $ u^{\alpha-\lambda}$, integrating over domain $U$, and using integration by parts we have
 \begin{align*}
\frac{\lambda}{\alpha}\frac{d }{dt}\int_U u^{\alpha} dx
&= \int_U \nabla\cdot \left(K(|\nabla u+Z(u)|)(\nabla u+Z(u)) \right) u^{\alpha-\lambda} dx +\int_U f u^{\alpha-\lambda}dx\\
&= - (\alpha-\lambda)\int_U K(|\nabla u+Z(u)|)(\nabla u+Z(u))\cdot \nabla u ( u^{\alpha-\lambda-1})  dx \\
&\quad +\int_{\Gamma} K(|\nabla u+Z(u)|)(\nabla u+Z(u))\cdot \vec \nu u^{\alpha-\lambda}   d\sigma
  +\int_U f u^{\alpha-\lambda}dx.
\end{align*}

For the first integral on the right-hand side, we write
\beqs
(\nabla u+Z(u))\cdot \nabla u=|\nabla u+Z(u)|^2-(\nabla u+Z(u))\cdot Z(u).
\eeqs
For the boundary integral, we use the boundary condition in \eqref{upb}.
These result in
\begin{align*}
&\frac{\lambda}{\alpha}\frac{d }{dt}\int_U u^{\alpha} dx
= - (\alpha-\lambda)\int_U K(|\nabla u+Z(u)|)|\nabla u+Z(u)|^2u^{\alpha-\lambda-1}  dx\\
& +(\alpha-\lambda)\int_U K(|\nabla u+Z(u)|)(\nabla u+Z(u))\cdot Z(u) u^{\alpha-\lambda-1} dx+\int_\Gamma u^{\alpha-\lambda} B d\sigma +\int_U f u^{\alpha-\lambda}dx.
\end{align*}

Using relations \eqref{Zprop},  \eqref{i:ineq1} and \eqref{i:ineq2} for the first two integrals on the right-hand side, as well as \eqref{Bprop}, \eqref{fprop} for the last two,  we have 
\begin{align*}
\frac{\lambda}{\alpha}\frac{d }{dt}\int_U u^{\alpha} dx
&\le - (\alpha-\lambda)d_3 \int_U (|\nabla u+Z(u)|^{2-a}-1)u^{\alpha-\lambda-1}  dx\\
& \quad +(\alpha-\lambda)d_2d_0 \int_U (|\nabla u+Z(u)|+1)^{1-a} u^{\alpha-\lambda-1+\ellz}dx\\
&\quad +\int_\Gamma u^{\alpha-\lambda}(\varphi_1+\varphi_2 u^\ellb)d\sigma+\int_U (f_1+f_2 u^\ellf) u^{\alpha-\lambda}dx.
\end{align*}
Using \eqref{ee2}, we estimate
\beqs
(|\nabla u+Z(u)|+1)^{1-a}\le |\nabla u+Z(u)|^{1-a}+1.
\eeqs
Hence, we obtain
\begin{align*}
\frac{\lambda}{\alpha}\frac{d }{dt}\int_U u^{\alpha} dx
&\le - (\alpha-\lambda)d_3I_1+(\alpha-\lambda)d_3  I_2+(\alpha-\lambda)d_2d_0I_3+(\alpha-\lambda)d_2d_0I_4+I_5+I_6.
\end{align*}
where
\begin{align*}
I_1&= \int_U |\nabla u+Z(u)|^{2-a}u^{\alpha-\lambda-1}  dx,
& I_2&= \int_Uu^{\alpha-\lambda-1}  dx,\\
I_3&= \int_U |\nabla u+Z(u)|^{1-a} u^{\alpha-\lambda-1+\ellz}dx,
& I_4&= \int_U u^{\alpha-\lambda-1+\ellz} dx,\\
I_5&= \int_\Gamma (u^{\alpha-\lambda} \varphi_1+u^{\alpha-\lambda+\ellb} \varphi_2) d\sigma ,
& I_6&=\int_U (u^{\alpha-\lambda} f_1+u^{\alpha-\lambda+\ellf} f_2) dx .
\end{align*}

In calculations below, the positive constants $\bar C$, $C$, $C_\varep$, $C'_\varep$ are generic, with $\bar C$ independent of $\alpha$, while $C$ depending on $\alpha$, and $C_\varep$, $C'_\varep$ depending on $\alpha$ and $\varep$.

$\bullet$ Applying Young's inequality to the integrand of $I_3$ with powers $\frac{2-a}{1-a}$ and $2-a$ gives 
\beq\label{Ithree}
I_3\le \frac{d_3}{2d_2d_0} \int_U |\nabla u+Z(u)|^{2-a}u^{\alpha-\lambda-1}  dx +\bar C\int_U u^{\ellz(2-a)+\alpha-\lambda-1} dx
=\frac{d_3}{2d_2d_0} I_1+\bar C I_7,
\eeq
where %$\bar C_1>0$ is independent of $\alpha$, and 
\beq\label{k1}
I_7= \int_U u^{\mu_1} dx\quad\text{with } \mu_1=\ellz(2-a)+\alpha-\lambda-1.
\eeq

Note also, by \eqref{ee4}, that
$$ I_4\le I_2+I_7.
$$

Then 
\beq\label{st1}
\frac{\lambda}{\alpha}\frac{d }{dt}\int_U u^{\alpha} dx\le - (\alpha-\lambda)\frac{d_3}2 I_1+\bar C\alpha(I_2+I_7)+I_5+I_6.
\eeq
%where $ \bar C_2>0$ is independent of  $\alpha$.

$\bullet$ To estimate $I_1$ we use inequality 
$$|x|^{2-a} \le 2^{1-a}(|x-y|^{2-a} +|y|^{2-a} )\text{ which gives  }
 |x-y|^{2-a}\ge  2^{a-1} |x|^{2-a} - |y|^{2-a}\quad\forall x,y\in\R^n.$$
  Together with \eqref{Zprop}, it gives 
 \beq\label{Ione}
\begin{aligned}
I_1
&\ge \int_U (2^{a-1} |\nabla u|^{2-a}u^{\alpha-\lambda-1}-|Z(u)|^{2-a}u^{\alpha-\lambda-1} ) dx\\
&\ge 2^{a-1}\int_U |\nabla u|^{2-a}u^{\alpha-\lambda-1}dx-d_0^{2-a} \int_U u^{\alpha-\lambda-1+\ellz(2-a)}  dx\\
&\ge \frac12\int_U |\nabla u|^{2-a}u^{\alpha-\lambda-1}dx-d_0^{2-a} I_7.
\end{aligned}
\eeq
This and \eqref{st1} imply 
\beq\label{Imain}
\frac{\lambda}{\alpha}\frac{d }{dt}\int_U u^{\alpha} dx + (\alpha-\lambda)\frac{d_3}{4} \int_U |\nabla u|^{2-a}u^{\alpha-\lambda-1}dx
\le \bar C\alpha(I_2+I_7)+I_5+I_6.
\eeq
%where $ \bar C_3>0$ is independent of  $\alpha$.

$\bullet$ For $I_5$, using Young's  inequality, we have 
\beq\label{I5y}
I_5\le \int_\Gamma u^{\mu_2}d\sigma +\int_\Gamma u^{\mu_3}d\sigma +\|\varphi_1\|_{L^{q_1}(\Gamma)}^{q_1}
+ \|\varphi_2\|_{L^{q_2}(\Gamma)}^{q_2},
\eeq
where
 \beq\label{k2}
 \mu_2=p_1(\alpha-\lambda)\quad\text{and}\quad \mu_3=p_2(\alpha-\lambda+\ellb).
 \eeq
 %Under \eqref{alal}, we see that $\mu_2,\mu_3>\alpha\ge 2>2-\delta$.  

For the first two integrals on the right-hand side, applying \eqref{t7}  to $\alpha+r=\mu_2$ and $\alpha+r=\mu_3$ gives
\beq\label{b1}
\begin{aligned}
\int_\Gamma u^{\mu_2}d\sigma +\int_\Gamma u^{\mu_3}d\sigma 
&\le \varep \int_U |u|^{\alpha-\lambda-1}|\nabla u|^{2-a} dx\\
 &\quad + c_1\int_U (u^{\mu_2}+u^{\mu_3})dx + C_\varep\int_U (u^{\mu_4}+ u^{\mu_5})dx,
\end{aligned}
\eeq
where $\varep>0$ is arbitrary,
\beq\label{mu1}
\begin{split}
\mu_4&=\alpha+\frac{(2-a)r+a-\delta}{1-a}=\alpha+\frac{(2-a)((p_1-1)\alpha-p_1\lambda)+a-\delta}{1-a}>\alpha,\\
\mu_5&=\alpha+\frac{(2-a)r+a-\delta}{1-a}=\alpha+\frac{(2-a)((p_2-1)\alpha+p_2(-\lambda+\ellb))+a-\delta}{1-a}>\alpha.
\end{split}
\eeq
(The validity of \eqref{b1} will be verified later.)
%Since $ \nu_1<\mu_i$ \mu_3<\alpha<\max\{\mu_4, \mu_5\}$, applying inequality \eqref{ee4} to the right hand side of \eqref{b1} yields  
Then
\beq\label{I5est}
\begin{aligned}
 I_5
&\le  
\varep \int_U |u|^{\alpha-\lambda-1}|\nabla u|^{2-a} dx  + C_\varep\int_U  \sum_{i=2}^5 u^{\mu_i}dx + \|\varphi_1\|_{L^{q_1}(\Gamma)}^{q_1} +\|\varphi_2\|_{L^{q_2}(\Gamma)}^{q_2}.
\end{aligned}
\eeq

$\bullet$ For $I_6$, using Young's inequalities we have  
\beq\label{I6est}
I_6\le  \int_U (u^{\mu_6}+ u^{\mu_7})dx + \|f_1\|_{L^{q_3}}^{q_3}
+  \|f_2\|_{L^{q_4}}^{q_4},
\eeq
where 
\beq\label{k4}
\mu_6=p_3(\alpha-\lambda)\quad\text{and}\quad
\mu_7= p_4(\alpha-\lambda+\ellf).\eeq

Note for $1\le i\le 7$ that $\nu_1\le \mu_i\le \nu_3$, where, referring to \eqref{k1}, \eqref{mu1} and \eqref{k4}, 
\beq\label{mu3}
\nu_3=\max \{  \mu_i:1\le i\le 7 \}.
\eeq
Hence by using \eqref{ee4}, $u^{\mu_i}\le u^{\nu_1}+u^{\nu_3}$  for all $1\le i \le 7$.
Combining this with \eqref{Imain}, \eqref{I5est}, and \eqref{I6est} yields \beq\label{Imain2}
\frac{\lambda}{\alpha}\frac{d }{dt}\int_U u^{\alpha} dx
+\Big[(\alpha-\lambda)\frac {d_3}4 -\varep\Big] \int_U |\nabla u|^{2-a}u^{\alpha-\lambda-1}dx
\le C_\varep \int_U (u^{\nu_1}+ u^{\nu_3})dx +\Upsilon(t),
\eeq

Using H\"{o}lder's inequality we have 
\beq\label{Itwo}
\int_U u^{\nu_1}  dx \le C \left( \int_U u^{\alpha}  dx\right)^\frac{\nu_1}{\alpha}= C \|u\|_{L^{\alpha}}^{\nu_1}.
\eeq

Applying \eqref{t6} to $\alpha+r=\nu_3$  gives
\beq\label{Imu3}
C_\varep\int_U u^{\nu_3}dx \le \varepsilon \int_U |u|^{\alpha-\lambda-1}|\nabla u|^{2-a} dx
+ C'_\varep\|u\|_{L^\alpha}^{\nu_3}+C'_\varep\|u\|_{L^\alpha}^{\nu_2}.
\eeq
(Again, we will verify the validity of \eqref{Imu3} later.)
Combining \eqref{Imain2} with \eqref{Itwo} and \eqref{Imu3}, also using inequality \eqref{ee4} with powers $\nu_1\le \nu_3\le \nu_2$ in dealing with $\|u\|_{L^\alpha}$, we have 
\beq\label{Imain3}
\frac{\lambda}{\alpha}\frac{d }{dt}\int_U u^{\alpha} dx
+\Big[(\alpha-\lambda)\frac {d_3}4 -2\varep\Big] \int_U |\nabla u|^{2-a}u^{\alpha-\lambda-1}dx
\le C'_\varep \|u\|_{L^{\alpha}}^{\nu_1} +C'_\varep \|u\|_{L^{\alpha}}^{\nu_2} +\Upsilon(t).
\eeq

Choosing $\varepsilon$ sufficiently small in \eqref{Imain2}, we obtain \eqref{est10}.

It remains to check conditions for inequalities \eqref{I5est} and \eqref{Imu3} to hold.
Inequality \eqref{b1} is valid under conditions 
 $ \mu_2, \mu_3>\alpha$ and $\alpha >2-\delta$, which, thanks to $\alpha\ge 2$, is equivalent to 
  \beq\label{Alplarge1}
 \alpha> \max\{q_1\lambda, q_2(\lambda -\ellb) \}.
 \eeq

The conditions for \eqref{Imu3} are $\nu_3>\alpha$, $\alpha\ge 2-a$, and 
\beq \label{Alplarge2}
\alpha>\frac{n(r+a-\delta)}{2-a}=\frac{n(\nu_3-\alpha+a-\delta)}{2-a}.
\eeq

Based on definition \eqref{mu3} of $\nu_3$, we consider the following cases.
 
 (i)  If $\nu_3=\mu_1$, condition \eqref{Alplarge2} becomes 
\beq\label{p11}
\alpha>n(\ellz-1).
\eeq

(ii)  If $\nu_3=\mu_2$, condition \eqref{Alplarge2} becomes 
\beq\label{p22}
\alpha>\frac{-p_1\lambda+a-\delta}{\kappa_f-p_1}.
\eeq

 (iii)  If $\nu_3=\mu_3$, condition \eqref{Alplarge2} becomes 
%$$\alpha>\frac{n(p_2(\alpha-\lambda+\ellb)-\alpha+a-\delta)}{2-a},$$
%which is equivalent to
\beq\label{p33}
\alpha>\frac{p_2(-\lambda+\ellb)+a-\delta}{\kappa_f-p_2}.
\eeq

(iv) If $\nu_3=\mu_4$,   condition \eqref{Alplarge2} becomes 
\beq\label{p44}
\alpha>\frac{-p_1\lambda+(a-\delta)}{\kappa_B-p_1}.
\eeq

(v) If $\nu_3=\mu_5$, condition \eqref{Alplarge2} becomes 
\begin{align*}
\alpha
&>\frac{n}{2-a}\Big(\frac{(2-a)((p_2-1)\alpha+p_2(-\lambda+\ellb))+a-\delta}{1-a}+a-\delta\Big)\\
&=\frac{n}{1-a}\Big((p_2-1)\alpha+p_2(-\lambda+\ellb)+a-\delta\Big),
\end{align*}
thus,
\beq\label{p55}
\alpha>\frac{p_2(-\lambda+\ellb)
+a-\delta}{\kappa_B-p_2}.
\eeq

(vi) If $\nu_3=\mu_6$, condition \eqref{Alplarge2} becomes  
\beq\label{p66}
\alpha>\frac{-p_3\lambda+a-\delta}{\kappa_f-p_3}.
\eeq

(vii)  If $\nu_3=\mu_7$, condition \eqref{Alplarge2} becomes 
%$$\alpha>\frac{n(p_4(\alpha-\lambda+\ellf)-\alpha+a-\delta)}{2-a},$$
%which is equivalent to
\beq\label{p77}
\alpha>\frac{p_4(-\lambda+\ellf)+a-\delta}{\kappa_f-p_4}.
\eeq

In summary, the conditions \eqref{Alplarge1}  and  \eqref{p11}--\eqref{p77}  are equivalent to $\alpha>\eta_0$, which is our choice of $\alpha$. The proof is complete.
\end{proof}
%==========================================================%

We then obtain the estimates in terms of initial and boundary data  in the next theorem.

%==========================================================%
\begin{theorem}\label{lathm}
Let $\alpha$, $C_{0}$, $\theta$   be as in Lemma \ref{lalem} and $\Upsilon(t)$ be defined by \eqref{Upsilondef}.
 If $T>0$ satisfies
\beq\label{Tcond1}
\int_0^T (1+\Upsilon(t))dt<C_{1}\cdot\Big(1+\int_U u_0(x)^\alpha dx\Big)^{-\nu_4},
\eeq
where $$\nu_4=\nu_4(\alpha)\eqdef \frac{(2-a)\theta}{n(1-\theta)}\quad\text{and}\quad C_{1}=C_{1}(\alpha)\eqdef \frac{1}{4C_{0}\nu_4},$$
then for all $t\in[0,T]$
\beq\label{La}
\int_U u^\alpha(x,t)dx \le  \Big\{\Big(1+\int_Uu_0(x)^\alpha dx\Big)^{-\nu_4}- \frac 1{C_{1}} \int_0^t (1+\Upsilon(\tau)) d\tau \Big\}^{-\frac{1}{\nu_4}}.
\eeq 

In particular, if $T>0$ satisfies 
\beq\label{Tcond2}
\int_0^T (1+\Upsilon(t))dt\le C_{1}\cdot(1-2^{-\nu_4}) \Big(1+\int_U u_0(x)^\alpha dx\Big)^{-\nu_4},
\eeq
then
\beq\label{La2}
\int_U u^\alpha(x,t)dx \le 2\Big(1+\int_U u_0^\alpha(x)dx\Big)\quad\text{for all}\quad t\in[0,T],
\eeq 
and
%\item[\rm (iii)] Also, (want this?)
\beq\label{La3}
\int_0^T \int_U |\nabla u|^{2-a}u^{\alpha-\lambda-1} dxdt \le  2(1+1/\nu_4)\Big(1+\int_U u_0(x)^{\alpha} dx\Big).
\eeq
\end{theorem}
%==============================%
\begin{proof}
In  \eqref{est10}, by  applying Young's inequality to the right-hand side we find that 
\beq
\begin{aligned}
\label{DIneq}
\frac{d }{dt}\int_U u^{\alpha} dx  +\int_U |\nabla u|^{2-a}u^{\alpha-\lambda-1}dx
&\le 2C_{0} \cdot(1+\|u(t)\|_{L^\alpha}^{\nu_2} +\Upsilon(t))\\
&\le 4C_{0}\cdot(1+\Upsilon(t)) \Big(1+\int_Uu^\alpha dx\Big)^{\frac{\nu_2}{\alpha}}.
\end{aligned}
\eeq
Let $V(t)=1+\int_Uu^\alpha(x,t) dx.$
It follows from \eqref{DIneq} that 
\beqs%\label{322}
V'(t)\le 4C_{0}\cdot(1+\Upsilon(t))V(t)^{\frac{\nu_2}{\alpha}}=4C_{0}\cdot(1+\Upsilon(t))V(t)^{\nu_4+1}.
\eeqs
Solving this differential inequality, under condition \eqref{Tcond1} yields the estimate \eqref{La} for all $t\in[0,T]$.
%gives 
%\beqs\label{V}
%V(t)\le \Big\{V(0)^{-\nu_4}- 4C_{0}\nu_4\int_0^t (1+\Upsilon(\tau)) d\tau \Big\}^{-\frac{1}{\nu_4}},
%\eeqs
%for all $t\in[0,T]$, with $T>0$ satisfying \eqref{Tcond1}. Hence we obtain estimate \eqref{La}.

Now, let $T>0$ satisfy \eqref{Tcond2}. Then we have from \eqref{La} that
\beq\label{V2}
V(t)\le 2V(0)\quad \forall t\in[0,T],
\eeq
and hence, estimate \eqref{La2} follows.

 Next, integrating \eqref{DIneq} in time from $0$ to $t$ and using \eqref{V2},  \eqref{Tcond2} we have
\begin{align}
& \int_0^T \int_U |\nabla u|^{2-a}u^{\alpha-\lambda-1}dt \le  \int_U u_0(x)^{\alpha} dx + 4C_{0}\int_0^T \Big(1+\int_U u^{\alpha} dx\Big)^{\nu_4+1} (1+\Upsilon(t)) dt \notag\\
& \le  \int_U u_0(x)^{\alpha} dx + 4C_{0}\cdot\Big(1+\int_U u_0(x)^{\alpha} dx\Big)^{\nu_4+1}\int_0^T  (1+\Upsilon(t)) dt \notag\\
&\le  \int_U u_0(x)^{\alpha} dx + 4C_{0}C_{1}\cdot(1-2^{-\nu_4})\Big(1+\int_U u_0(x)^{\alpha} dx\Big). \label{xtgrad}
\end{align}
Note that
$$
4C_{0}C_{1}\cdot(1-2^{-\nu_4})=\frac{1-2^{-\nu_4}}{\nu_4}\le \frac{1}{\nu_4}.
$$
Thus we obtain \eqref{La3} from \eqref{xtgrad}.
\end{proof}

%%%%%%%%%%%%%%%%%%%%%%%%%%%%%%%%%%%%%%%%%%%%%%%%%%%%%%%%%%%%%%%%%%%%%
\myclearpage
\section{Maximum estimates}\label{maxsec}
%%%%%%%%%%%%%%%%%%%%%%%%%%%%%%%%%%%%%%%%%%%%%%%%%%%%%%%%%%%%%%%%%%%%%

In this section, we use Moser's iteration to estimate the $L^\infty$-norm of the solution $u(x,t)$. 

For $i=1,2,3,4$, let $q_i$ and $p_i$ be as in section \ref{lasec}.
Assume \eqref{pscond}. Let  
\beq\label{fixedpower}
\kappa_*=\max\Big\{ \frac{(2-a)p_1-1}{1-a},\frac{(2-a)p_2-1}{1-a},p_3,p_4 \Big\}, 
\eeq
\beqs
\eta_1= \max\Big\{\frac{p_1\lambda-a+\delta}{(2-a)p_1-1}, \frac{p_2(\lambda-\ellb)-a+\delta}{(2-a)p_2-1}\Big\},
\eeqs
\beqs
\eta_2=\max\Big\{\ellz(2-a),\ellf,\frac{a-\delta+(2-a)p_2 \ellb}{1-a} \Big\}.
\eeqs

%================================
\begin{lemma} \label{caccio}
Given  $\tilde \kappa>\kappa_*$,  suppose
\beq\label{alp-Large}
\alpha > \max\{2, \alpha_*,\eta_1\} \quad\text{and}\quad \alpha \ge \frac{\eta_2}{\tilde\kappa-\kappa_*}.
\eeq

If $T>T_2>T_1\ge 0$ then 
\beq\label{Sest2}
\begin{aligned}
&\sup_{t\in[T_2,T]} \int_U u^\alpha(x,t)dx +\int_{T_2}^T\int_U |\nabla u(x,t)|^{2-a}u(x,t)^{\alpha-\lambda-1}dx dt\\
&\le  c_6(1+ T)\Big(1+\frac1{T_2-T_1}\Big)\alpha^2 \mathcal M_0\Big(\norm {u}_{L^{\tilde \kappa\alpha}(U\times(T_1,T))}^{\nu_5}+\norm {u}_{L^{\tilde \kappa\alpha}(U\times(T_1,T))}^{\nu_6}\Big),
\end{aligned}
\eeq
where $c_6\ge 1$ is independent of $\alpha$, $\tilde \kappa$, $T$, $T_1$ and  $T_2$, % the quantity $\mathcal M_0$ is defined by
\beq\label{Mbar}
\mathcal M_0=1+ \norm{\varphi_1}_{L^{q_1}(\Gamma_T)}^\frac{2-a}{1-a}
+ \norm{\varphi_2}_{L^{q_2}(\Gamma_T)}^\frac{2-a}{1-a}+\norm{f_1}_{L^{q_3}(Q_T)}+\norm{f_2}_{L^{q_4}(Q_T)},
\eeq
the positive powers $\nu_5$ and $\nu_6$ are defined by
\beq\label{mu12}
 \nu_5=\alpha -h_1 \quad\text{and}\quad \nu_6=\alpha +h_2
\eeq
with 
\beq\label{h12def}
\begin{aligned}
h_1&=\max\Big\{\lambda+1,\frac{p_1\lambda -a+\delta}{p_1(2-a)-1},\frac{p_2(\lambda-\ellb)-a+\delta}{p_2(2-a)-1}\Big\}>0,\\
h_2&=\max\Big\{ 0,\ellz(2-a)-\lambda-1,\ellf-\lambda, \ellb-\lambda, \frac{a-\delta-p_1\lambda}{p_1(2-a)-1},\frac{a-\delta-p_2(\lambda-\ellb)}{p_2(2-a)-1}\Big\}\ge 0.
\end{aligned}
\eeq
\end{lemma}
%===========================================%
\begin{proof}
Denote $Q_T =U\times[0,T]$ and $\Gamma_T = \Gamma\times [0,T]$.  Let $\xi=\xi(t)$ be a $C^1$-function  on $[0,T]$ with $\xi(0)=0$ and $0\le \xi(t)\le 1$. In calculations below, the generic positive constant $\bar C$ is independent of $\alpha$, $\tilde \kappa$, $T$, $T_1$ and  $T_2$.

Multiply the PDE in \eqref{upb} by test $u^{\alpha-\lambda}\xi^2$ and  integrating the resulting equation over $U$ give
 \begin{align*}
&\frac{\lambda}{\alpha}\frac{d }{dt}\int_U u^{\alpha} \xi^2 dx
-\frac{\lambda}{\alpha}\int_U 2u^{\alpha} \xi \xi' dx\\
&= \int_U \nabla\cdot \Big(K(|\nabla u+Z(u)|)(\nabla u+Z(u)) \Big) u^{\alpha-\lambda}\xi^2 dx +\int_U f u^{\alpha-\lambda } \xi^2 dx.
\end{align*}

Since $\xi(t)$ is independent of $x$, same estimates as in section \ref{lasec} give the following version of \eqref{Imain}
 \begin{align*}
&\frac{\lambda}{\alpha}\frac{d }{dt}\int_U u^{\alpha}\xi^2 dx
+\frac{(\alpha-\lambda)d_3}{4} \int_U |\nabla u|^{2-a}u^{\alpha-\lambda-1}\xi^2 dx 
\le C\alpha\int_U(u^{\nu_1}+u^{\mu_1}) \xi^2 dx\\
& +\int_\Gamma u^{\alpha-\lambda}(\varphi_1+\varphi_2 u^\ellb)\xi^2 d\sigma+\int_U (f_1+f_2 u^\ellf) u^{\alpha-\lambda}\xi^2 dx+\frac{\lambda}{\alpha}\int_U 2u^{\alpha} \xi \xi' dx.
\end{align*}
%where $C>0$ is independent of $\alpha$.

Integrating the previous inequality in time and applying H\"older's inequality give
\begin{multline}\label{Imain4}
\frac{\lambda}{\alpha}\sup_{t\in[0,T]}\int_U u^{\alpha}(x,t)\xi^2(t) dx
+\frac{(\alpha-\lambda)d_3}4 \iint_{Q_T} |\nabla u|^{2-a}u^{\alpha-\lambda-1} \xi^2 dxdt \\
\le C\alpha\iint_{Q_T} (u^{\nu_1}+u^{\mu_1})\xi^2  dxdt
 +\frac{\lambda}{\alpha}\iint_{Q_T} 2u^{\alpha} \xi \xi' dxdt +\sum_{i=1}^4 I_i 
\end{multline}
where 
\begin{align*}
I_i& =\Big (\iint_{\Gamma_T}  u^{\alpha_i}\xi^2d\sigma dt\Big)^{1/p_i} F_i    \quad\text{ with }
	F_i=\Big(\iint_{\Gamma_T} \varphi_i^{q_i}\xi^2 d\sigma dt\Big)^{1/q_i}  \text{ for } i=1,2, \\
%I_2& =\Big (\int_{\Gamma_T}  u^{\alpha_2}\xi^2d\sigma dt\Big)^{1/p_2} F_2\text{ with }
%	F_2=\Big(\int_{\Gamma_T} \varphi_2^{q_2}\xi^2 d\sigma dt\Big)^{1/q_2}, \\
I_3&=\Big(\iint_{Q_T} u^{\alpha_3}\xi^2 dxdt\Big)^{1/p_3} F_3 \quad\text{ with }
	F_3=\Big(\iint_{Q_T} f_1^{q_3}\xi^2dxdt\Big)^{1/q_3},\\
I_4&=\Big(\iint_{Q_T} u^{\alpha_4}\xi^2 dxdt\Big)^{1/p_4}F_4 \quad\text{ with }
	F_4= \Big(\iint_{Q_T} f_2^{q_4}\xi^2dxdt\Big)^{1/q_4},
\end{align*}
with
\begin{align*}
 \alpha_1&=p_1(\alpha-\lambda)>p_1>1,& \alpha_2&=p_2(\alpha-\lambda+\ellb)>p_2(1+\ellb)>1,\\
\alpha_3&=p_3(\alpha-\lambda)>p_3>1, & \alpha_4&=p_4(\alpha-\lambda+\ellf)>p_4(1+\ellf)>1.
\end{align*}

For $I_1$, by using the trace theorem \eqref{W11}, we have 
\beqs
 \iint_{\Gamma_T } u^{\alpha_1}\xi^2 d\sigma dt
\le c_1 \iint_{Q_T} u^{\alpha_1}\xi^2 dxdt +c_2\alpha_1 \iint_{Q_T} u^{\alpha_1-1}|\nabla u|\xi^2 dxdt.
\eeqs

In the last integral, writing $u^{\alpha_1-1}|\nabla u|$ as a product of $u^\frac{\alpha-\lambda-1}{2-a}|\nabla u|$ and $u^{\alpha_1-\frac{\alpha-a+\delta}{2-a}}$, and applying H\"older's inequality with powers $2-a$ and $(2-a)/(1-a)$, we obtain
\begin{multline}\label{i1x}
 \iint_{\Gamma_T } u^{\alpha_1}\xi^2 d\sigma dt
\le c_1 \iint_{Q_T} u^{\alpha_1}\xi^2 dxdt \\+c_2\alpha_1\Big (\iint_{Q_T} u^{\alpha-\lambda-1}|\nabla u|^{2-a}\xi^2  dxdt\Big)^\frac1{2-a}\Big(\iint_{Q_T} u^{m_1}\xi^2 dxdt\Big)^\frac{1-a}{2-a},
\end{multline}
where
\beqs
m_1=\frac{(2-a)\alpha_1-\alpha+a-\delta}{1-a}=\frac{((2-a)p_1-1)\alpha-p_1\lambda+a-\delta}{1-a}>0.
\eeqs

Denote $X(\alpha)=\iint_{Q_T} u^{\alpha}(x,t)\xi(t)  dx dt$. We will use the fact $0\le \xi^2\le \xi\le 1$.
By \eqref{i1x} and Young's inequality
\begin{align}
I_1&\le c_1^{1/p_1}X(\alpha_1)^{1/p_1} F_1 +(c_2\alpha_1)^{1/p_1} \Big(\iint_{Q_T} u^{\alpha-\lambda-1}|\nabla u|^{2-a}dxdt\Big)^\frac1{p_1(2-a)}X(m_1)^\frac{1-a}{p_1(2-a)}F_1 \notag\\
&\le c_1^{1/p_1}X(\alpha_1)^{1/p_1}F_1
 + \varep \iint_{Q_T} u^{\alpha-\lambda-1}|\nabla u|^{2-a}dxdt \notag\\
&\quad +[\varep^{-1}(c_2\alpha_1)^{2-a} ]^{\frac1{p_1(2-a)-1} } X(m_1)^\frac{1-a}{p_1(2-a)-1}F_1^\frac{p_1(2-a)}{p_1(2-a)-1}. \label{II1}
\end{align}

Similarly, 
\beq\label{II2}
\begin{aligned}
I_2&\le c_3X(\alpha_2)^{1/p_2}F_2 + \varep \iint_{Q_T} u^{\alpha-\lambda-1}|\nabla u|^{2-a}dxdt\\
&\quad + [\varep^{-1}(c_4\alpha_2)^{2-a} ]^{\frac1{p_2(2-a)-1} } X(m_2)^\frac{1-a}{p_2(2-a)-1}F_2^\frac{p_2(2-a)}{p_2(2-a)-1},
\end{aligned}
\eeq
where
\beqs
m_2=\frac{(2-a)\alpha_2 -\alpha+a-\delta}{1-a}=\frac{((2-a)p_2-1)\alpha+p_2(\ellb-\lambda)+a-\delta}{1-a}>0.
\eeqs

We now  choose $\xi(t)$ such that 
\beq\label{xiprop} 
\xi(t)=0 \text{ for }0\le t\le T_1, \ 
\xi(t)=1 \text{ for }T_2\le t\le T, \text{ and }
 0\le \xi'(t)\le \frac{2}{T_2-T_1} \text{ for }0\le t\le T.
\eeq 
Then
\beq\label{IIx}
\frac{\lambda}{\alpha}\iint_{Q_T} 2u^{\alpha} \xi \xi' dxdt\le \frac{\bar C}{T_2-T_1}X(\alpha).
\eeq

Therefore, with $\varep=(\alpha-\lambda)d_3/16,$ it follows \eqref{Imain4} and the above estimates \eqref{II1}, \eqref{II2}, \eqref{IIx} that 
\beq\label{Imain5}
\begin{aligned}
&\frac{\lambda}{\alpha}\sup_{t\in[0,T]}\int_U u^{\alpha}(x,t)\xi^2(t) dx
  +\frac{(\alpha-\lambda)d_3}8 \iint_{Q_T} |\nabla u|^{2-a}u^{\alpha-\lambda-1} \xi^2 dx dt 
  \le \bar C\alpha (X(\mu_1)+ X(\nu_1))\\
&+\frac{\bar C}{T_2-T_1}X(\alpha)+ \bar C\sum_{i=1}^4 X(\alpha_i)^{1/p_i}F_i
+\bar C \sum_{j=1}^2 \Big(\frac{\alpha_j^{2-a}}{\alpha-\lambda} \Big)^{\frac1{p_j(2-a)-1} } X(m_j)^\frac{1-a}{p_j(2-a)-1}F_j^\frac{p_j(2-a)}{p_j(2-a)-1}.
\end{aligned}
\eeq

We bound the left-hand side of \eqref{Imain5} from below by  
\beqs
I_0\eqdef \frac{\lambda}{\alpha}\sup_{t\in[0,T]}\int_U u^{\alpha}(x,t)\xi^2(t) dx
 +\frac{d_3}8 \iint_{Q_T} |\nabla u|^{2-a}u^{\alpha-\lambda-1} \xi^2 dx dt.
 \eeqs
 
We now bound the right-hand side of \eqref{Imain5}. We will estimate many $X(\cdot)$-terms by using
\beqs
J_0\eqdef \norm {u}_{L^{\tilde \kappa\alpha}(U\times(T_1,T))}.
\eeqs

\emph{Claim.} Under condition \eqref{alp-Large}, one has  
\beq\label{powers}
0<\mu_1,\alpha_1,\alpha_2,
     \alpha_3, \alpha_4,m_1,m_2 \le \tilde \kappa\alpha .
\eeq

We accept this claim at the moment. By H\"older's inequality if $0<\beta\le\tilde \kappa \alpha$ then
\beq\label{X}
X(\beta)\le \int_{T_1}^T \int_U u^\beta  dxdt \le J_0^\beta |Q_T|^{1-\frac{\beta}{\tilde \kappa \alpha}} \le J_0^\beta (1+|Q_T|).
\eeq

Note for  $j=1,2$ that 
\beq\label{pjs}
\frac{1-a}{p_j(2-a)-1} <1\quad\text{and}\quad \frac{p_j(2-a)}{p_j(2-a)-1}<\frac {2-a}{1-a}.
\eeq
Using also $\alpha-\lambda\ge \alpha/2$ and $\alpha_1,\alpha_2\le (p_1+p_2)(1+\ellb)\alpha$, then for $j=1,2$,
\beqs
 \Big(\frac{\alpha_j^{2-a}}{\alpha-\lambda} \Big)^{\frac 1{p_j(2-a)-1}}
 \le  \Big(2[(p_1+p_2)(1+\ellb)]^{2-a} \alpha^{1-a}\Big)^{\frac 1{p_j(2-a)-1}}
  %\le  C\alpha^{\frac{1-a}{p_j(2-a)-1} }
  \le  C\alpha.
\eeqs
%where  $ C>0$ is independent of $\alpha$, $\tilde\kappa$.

Thus, we have from \eqref{Imain5}, \eqref{X} that
\beqs
I_0\le \bar C(1+|Q_T|)\Big (\alpha J_0^{\mu_1}
+\alpha J_0^{\nu_1}
+ \frac{1}{T_2-T_1}  J_0^\alpha   +\sum_{i=1}^4 J_0^\frac{\alpha_i}{p_i} F_i
+\alpha \sum_{j=1}^2 J_0^\frac{m_j(1-a)}{p_j(2-a)-1}F_j^\frac{p_j(2-a)}{p_j(2-a)-1}\Big).
 \eeqs

Calculating the powers of $J_0$ gives
\beqs
\frac{\alpha_1}{p_1}=\frac{\alpha_3}{p_3}=\alpha-\lambda,\quad
\frac{\alpha_2}{p_2}=\alpha-\lambda+\ellb,\quad
\frac{\alpha_4}{p_4}=\alpha-\lambda+\ellf,
\eeqs
\beqs
\frac{m_1(1-a)}{p_1(2-a)-1}=\alpha+\frac{-p_1\lambda+a-\delta}{p_1(2-a)-1},\quad
\frac{m_2(1-a)}{p_2(2-a)-1}=\alpha+\frac{p_2(\ellb-\lambda)+a-\delta}{p_2(2-a)-1}.
\eeqs
Hence,
\beqs
I_0 \le \bar C(1+|Q_T|)\alpha\Big(1+\frac {1}{T_2-T_1}\Big)\tilde F\cdot\tilde J.
 \eeqs
 where
\begin{align*}
\tilde F &= 1+F_1 +F_2+F_1^\frac{p_1(2-a)}{p_1(2-a)-1}+ F_2^\frac{p_2(2-a)}{p_2(2-a)-1} +F_3+F_4, \\
\tilde J &=J_0^{\alpha +\ellz(2-a)-\lambda-1 }+ J_0^{\alpha-\lambda-1}+J_0^\alpha+J_0^{\alpha-\lambda}+J_0^{\alpha+\frac{-p_1\lambda+a-\delta}{p_1(2-a)-1}}+ J_0^{\alpha +\ellb-\lambda }\\
&\quad  +J_0^{\alpha+\frac{p_2(\ellb-\lambda)+a-\delta}{p_2(2-a)-1}}+J_0^{\alpha +\ellf-\lambda}.
\end{align*}
%with 
%$
%\mathcal{S}= \norm {u}_{L^{\tilde \kappa\alpha}(U\times(T_1,T))}.
%$
Using the second inequality in \eqref{pjs} and inequality   \eqref{ee5}, we simply estimate
\beqs
\tilde F \le 3(1+F_1^\frac{2-a}{1-a}+ F_2^\frac{2-a}{1-a} +F_3+F_4)
\le  3\mathcal M_0.
\eeqs
Also, by using \eqref{ee4},
\beqs
\tilde J \le 8 (J_0^{\alpha-h_1} +J_0^{\alpha+h_2}).
\eeqs
Therefore 
\beq\label{I0i}
I_0  \le \bar C(1+|Q_T|)\alpha\Big(1+\frac {1}{T_2-T_1}\Big) \mathcal M_0(J_0^{\alpha-h_1} +J_0^{\alpha+h_2}). 
\eeq
By \eqref{I0i} and definition of $I_0$ we have
\beqs
\begin{split}
&\sup_{t\in[0,T]}\int_U u^{\alpha}(x,t)\xi^2(t) dx
 + \iint_{Q_T} |\nabla u|^{2-a}u^{\alpha+\lambda-1} \xi^2 dx dt\\
& \le \bar C(1+|Q_T|)\alpha^2\Big(1+\frac {1}{T_2-T_1}\Big) \mathcal M_0(J_0^{1-h_1/\alpha} +J_0^{1+h_2/\alpha}).
 \end{split}
\eeqs

This estimate, together with property \eqref{xiprop} of function $\xi(t)$, and
the fact $1+|Q_T|\le (1+|U|)(1+T)$,  implies \eqref{Sest2}.

Finally, we verify the claim \eqref{powers}. Dividing \eqref{powers} by $\alpha$ gives an equivalent statement
\begin{multline}\label{ktil2}
       1+  \frac{\ellz(2-a)-\lambda-1}{\alpha},\
     p_1-\frac{\lambda}{\alpha},\
     p_2+\frac{\ellb-\lambda}{\alpha},\
     p_3-\frac{\lambda}{\alpha},\
     p_4+\frac{\ellf-\lambda}{\alpha},\\
     \frac{(2-a)p_1-1}{1-a} +\frac{-(2-a)p_1\lambda+a-\delta}{\alpha(1-a)},\
     \frac{(2-a)p_2-1}{1-a}+ \frac{(2-a)p_2(\ellb-\lambda)+a-\delta}{\alpha(1-a)} 
 \le \tilde\kappa. 
\end{multline}
Note for $i=1,2$ that $p_i\le \frac{(2-a)p_i-1}{1-a}$. Then by definition \eqref{fixedpower} of $\kappa_*$,
\beq\label{ktilcond}
       1, p_1,p_2,p_3,p_4,\frac{(2-a)p_1-1}{1-a},\frac{(2-a)p_2-1}{1-a}\le \kappa_*. 
\eeq 
Hence, we can bound the left-hand side of \eqref{ktil2} by
\begin{align}
&\kappa_*+\max\Big\{ \frac{\ellz(2-a)-\lambda-1}{\alpha},\frac{\ellb-\lambda}{\alpha},
        \frac{\ellf-\lambda}{\alpha},
        \frac{a-\delta-(2-a)p_1\lambda}{\alpha(1-a)}, 
        \frac{a-\delta+(2-a)p_2(\ellb-\lambda)}{\alpha(1-a)} \Big\} \notag\\
&\le \kappa_*+\frac{\eta_2}\alpha.\label{conds}
\end{align}
Since $\tilde\kappa>\kappa_*$, we choose a sufficient condition for \eqref{ktil2} to be
$
\eta_2/\alpha\le \tilde\kappa-\kappa_*,
$
which is satisfied by \eqref{alp-Large}. The proof is complete.
\end{proof}
%================================================================%

%================================================================%
\begin{lemma}\label{GLk}
Let $\tilde \kappa$, $\alpha$,   $\nu_5$, $\nu_6$  and $\mathcal M_0$ be as in Lemma \ref{caccio}. 
If $T>T_2>T_1\ge 0$  then
%\Enote{erase this: and $\alpha>0$} 
\beq\label{bfinest2}
\| u\|_{L^{\kappa\alpha}(U\times (T_2,T))}
\le A_\alpha^\frac 1{\alpha}\Big( \| u\|_{L^{\tilde \kappa \alpha}(U\times(T_1,T))}^{\nu_5}+\| u\|_{L^{\tilde\kappa \alpha}(U\times(T_1,T))}^{\nu_7}\Big)^\frac 1{\alpha},
\eeq
where $\kappa$ is defined by \eqref{kappadef},  
%$\nu_5=\nu_5(\alpha)$ is defined in \eqref{mu12},
\beq\label{tilrsdef}
%\tilde r=\tilde r(\alpha)\eqdef \nu_5,\quad 
\nu_7=\nu_7(\alpha)\eqdef  \frac{\alpha\nu_6}{\alpha+\delta-a},
\eeq
and 
\beq\label{Atildef2}
A_\alpha = c_7 (1+ T)^2 \Big(1+\frac1{T_2-T_1}\Big)^2 \alpha^{6-a} \mathcal M_0^2
\eeq
with  $c_7\ge 1$ independent of $\alpha$, $\tilde \kappa$, $T$, $T_1$ and $T_2$.
\end{lemma}
%================================================================%
\begin{proof} We follow the proof of Proposition 6.2 in \cite{CHK1}. 
We use Sobolev inequality  \eqref{parabwki} in Lemma \ref{Parabwk}:
\beq\label{beginest2}
\begin{aligned}
%J \eqdef  \Big(\int_{T_2}^T \int_U |u|^{ \kappa\alpha}dx dt \Big)^{\frac 1{\kappa\alpha}}
\| u\|_{L^{\kappa\alpha}(U\times (T_2,T))}
&\le \hat C^\frac{1}{\kappa\alpha} \Big\{I^{\tilde \theta}\cdot \sup_{t\in[T_2,T]}\Big( \int_U|u(x,t)|^{\alpha}dx \Big)^{1-\tilde \theta}\Big\}^{\frac 1\alpha},
\end{aligned}
\eeq
where $\hat C=c_5\alpha^{2-a}$ with $c_5\ge 1$, the numbers $\tilde \theta$ and $\kappa$ are defined in \eqref{kappadef} and 
\beqs
\begin{aligned}
I=\Big[  \int_{T_2}^T\int_U |u|^{\alpha+\delta-a}dxdt+\int_{T_2}^T\int_U |u|^{\alpha-\lambda-1}|\nabla u|^{2-a}dxdt \Big]^{\frac {\alpha}{\alpha+\delta-a}} 
\end{aligned}
\eeqs
Note that 
\beq\label{C5}
\frac {\alpha}{\alpha+\delta-a}\le 2 \quad \text{and}\quad 2^{\frac {\alpha}{\alpha+\delta-a}-1}\le 2.
\eeq
Then applying inequality \eqref{ee3}, we find that   
\beq\label{Iest2}
I\le 2 \Big(\int_{T_2}^T\int_U |u|^{\alpha+\delta-a}dxdt\Big)^{\frac {\alpha}{\alpha+\delta-a}}+2\Big(\int_{T_2}^T\int_U |u|^{\alpha-\lambda-1}|\nabla u|^{2-a}dxdt \Big)^{\frac {\alpha}{\alpha+\delta-a}}.
\eeq
Applying H\"older's inequality to the first integral on the right-hand side of \eqref{Iest2} with conjugate exponents $\frac{\tilde\kappa\alpha}{\alpha+\delta-a}$ and $\frac{\tilde\kappa\alpha}{(\tilde\kappa-1)\alpha+a-\delta}$, we get 
\beq\label{Iest3}
I \le 2 |Q_T|^{(1-\frac{\alpha+\delta-a}{\tilde\kappa  \alpha}) \frac {\alpha}{\alpha+\delta-a}}\Big(\int_{T_2}^T\int_U |u|^{\tilde\kappa \alpha}dx dt\Big)^\frac1{\tilde\kappa}+2\Big(\int_{T_2}^T\int_U |u|^{\alpha-\lambda-1}|\nabla u|^{2-a}dx dt \Big)^{\frac {\alpha}{\alpha+\delta-a}}.
\eeq

Next, we use \eqref{Sest2} to estimate the second term on the right-hand side of \eqref{Iest3}.
We denote
\beq\label{newMdef}
\mathcal{S}= \norm {u}_{L^{\tilde \kappa\alpha}(U\times(T_1,T))}
\quad \text{and}\quad \mathcal M=c_6(1+ T)\Big(1+\frac1{T_2-T_1}\Big)\alpha^2 \mathcal M_0,
\eeq
where $c_6$ is the positive constant in \eqref{Sest2}.
Then combining $\eqref{beginest2}$ and $\eqref{Iest3}$ yields 
\begin{align*}
%J
\| u\|_{L^{\kappa\alpha}(U\times (T_2,T))}
&\le  \hat C^\frac{1}{\kappa\alpha} \Big\{ \Big [ 2(1+ |Q_T|)^{ \frac {\alpha}{\alpha+\delta-a}}\mathcal{S}^\alpha+2(\mathcal M(\mathcal{S}^{\nu_5}+\mathcal{S}^{\nu_6}))^{\frac{\alpha}{\alpha+\delta-a}}\Big]^{\tilde \theta}(\mathcal M(\mathcal{S}^{\nu_5}+\mathcal{S}^{\nu_6}))^{1-\tilde \theta} \Big\}^{\frac 1{\alpha}}\\
&\le  \hat C^\frac{1}{\kappa\alpha} \Big\{  \Big [ 2 (1+|Q_T|)^2 \mathcal{S}^\alpha+2(\mathcal M(\mathcal{S}^{\nu_5}+\mathcal{S}^{\nu_6}))^{\frac{\alpha}{\alpha+\delta-a}}+\mathcal M(\mathcal{S}^{\nu_5}+\mathcal{S}^{\nu_6})\Big]^{\tilde \theta+1-\tilde \theta} \Big\}^{\frac 1{\alpha}}\\
&= \hat C^\frac{1}{\kappa\alpha} \Big\{   2 (1+|Q_T|)^2 \mathcal{S}^\alpha+2\mathcal M^{\frac{\alpha}{\alpha+\delta-a}}(\mathcal{S}^{\nu_5}+\mathcal{S}^{\nu_6})^{\frac{\alpha}{\alpha+\delta-a}} +\mathcal M(\mathcal{S}^{\nu_5}+\mathcal{S}^{\nu_6})\Big\}^{\frac 1{\alpha}}.
\end{align*}
Since $\nu_5<\alpha<\nu_6$, we use \eqref{ee4} to estimate
$\mathcal{S}^\alpha\le \mathcal{S}^{\nu_5}+ \mathcal{S}^{\nu_6}.$
Also, by \eqref{ee3} and \eqref{C5},
\beqs
(\mathcal{S}^{\nu_5}+\mathcal{S}^{\nu_6})^{\frac{\alpha}{\alpha+\delta-a}}\le 2(\mathcal{S}^{\frac{\nu_5\alpha}{\alpha+\delta-a}}+\mathcal{S}^{\frac{\nu_6\alpha}{\alpha+\delta-a}}) 
\eeqs
Thus, we have
\beqs
\| u\|_{L^{\kappa\alpha}(U\times (T_2,T))}
%\le \hat C^\frac{1}{\kappa\alpha} \Big\{  [ 2  (1+|Q_T|)^2  +\mathcal M](\mathcal{S}^{\nu_5}+\mathcal{S}^{\nu_6})
%+2\mathcal M^{\frac{\alpha}{\alpha+\delta-a}}(\mathcal{S}^{\nu_5}+\mathcal{S}^{\nu_6})^{\frac{\alpha}{\alpha+\delta-a}} \Big\}^{\frac 1{\alpha}}\\
\le  \hat C^\frac{1}{\kappa\alpha} \Big\{  \Big( 2  (1+|Q_T|)^2  +\mathcal M\Big)(\mathcal{S}^{\nu_5}+\mathcal{S}^{\nu_6})
+4\mathcal M^{\frac{\alpha}{\alpha+\delta-a}} (\mathcal{S}^{\frac{\nu_5\alpha}{\alpha+\delta-a}}+\mathcal{S}^{\frac{\nu_6\alpha}{\alpha+\delta-a}}) \Big\}^{\frac 1{\alpha}}.
\eeqs
Using $\mathcal{S}^{\nu_6},\mathcal{S}^{\frac{\nu_5\alpha}{\alpha+\delta-a}}\le \mathcal{S}^{\nu_5}+\mathcal{S}^{\frac{\nu_6\alpha}{\alpha+\delta-a}}$, we obtain
\beq\label{newJ}
\| u\|_{L^{\kappa\alpha}(U\times (T_2,T))}\le  \hat C^\frac{1}{\kappa\alpha} \Big\{ 8\Big(   (1+|Q_T|)^2  +\mathcal M+\mathcal M^{\frac{\alpha}{\alpha+\delta-a}}\Big)\Big(\mathcal{S}^{\nu_5}+\mathcal{S}^{\frac{\nu_6\alpha}{\alpha+\delta-a}}\Big) \Big\}^{\frac 1{\alpha}}.
\eeq

Note also that $\hat C^\frac1\kappa\le \hat C$,
then we obtain from \eqref{newJ} that 
\beq\label{uB}
\| u\|_{L^{\kappa\alpha}(U\times (T_2,T))}
\le \mathcal M_1^\frac 1{\alpha}\Big( \| u\|_{L^{\tilde \kappa \alpha}(U\times(T_1,T))}^{\nu_5}+\| u\|_{L^{\tilde\kappa \alpha}(U\times(T_1,T))}^{\nu_7}\Big)^\frac 1{\alpha},
\eeq
where
\beq\label{Atildef}
\mathcal M_1 = 8c_5\alpha^{2-a}  [  (1+|Q_T|)^2+\mathcal M+\mathcal M^{\frac{\alpha}{\alpha+\delta-a}}].\eeq

We estimate $\mathcal M_1$. Note  that  $ \mathcal M_0,\mathcal M,\mathcal M_1\ge1$.
Recall from \eqref{C5} that the last power of $\mathcal M$ in \eqref{Atildef} is less than or equal to $2$.
%Since  $\alpha\ge 2$, we have 
%$$\frac{\alpha}{\alpha+\delta-a}\le 2\quad\text{ and }\quad\frac{(\tilde\kappa-1)\alpha+a-\delta}{\tilde\kappa(\alpha+\delta-a)}
%\le \frac{2(\tilde\kappa-1)+a-\delta}{\tilde\kappa(2+\delta-a)}\le \frac{2\tilde\kappa}{\tilde \kappa}=2.$$
Then\beq\label{ti9}
\begin{aligned}
\mathcal M_1
&\le  8c_5\alpha^{2-a} \Big((1+|U|)^2(1+T)^2+2\mathcal M^2\Big)\\
&\le  8c_5\alpha^{2-a} \Big((1+|U|)^2(1+T)^2+2\Big\{c_6\alpha^2 \Big(1+\frac1{T_2-T_1}\Big)(1+ T)\mathcal M_0\Big\}^2\Big)\\
&\le 8c_5\Big ((1+|U|)^2+2c_6^2\Big) \alpha^{6-a}\Big(1+\frac1{T_2-T_1}\Big)^2(1+ T)^2 \mathcal M_0^2.
\end{aligned}
\eeq
Hence we obtain \eqref{bfinest2} from \eqref{uB} and \eqref{ti9}.
\end{proof}

%=======================================================%

Next, we apply Moser's iteration. For that we recall the following lemma on finding the upper bounds of certain numeric sequences.

%======================================================%
\begin{lemma}[\cite{CHK1}, Lemma A.2]\label{Genn}
Let $y_j\ge 0$, $\kappa_j>0$, $s_j \ge r_j>0$ and $\omega_j\ge 1$  for all $j\ge 0$.
Suppose there is $A\ge 1$ such that
\beqs
y_{j+1}\le A^\frac{\omega_j}{\kappa_j} (y_j^{r_j}+y_j^{s_j})^{\frac 1{\kappa_j}}\quad\forall j\ge 0.
\eeqs
Denote $\beta_j=r_j/\kappa_j$ and $\gamma_j=s_j/\kappa_j$.
Assume
\beqs
\bar\alpha \eqdef \sum_{j=0}^{\infty}  \frac{\omega_j}{\kappa_j}<\infty\text{ and the products }
\prod_{j=0}^\infty \beta_j, 
\prod_{j=0}^\infty \gamma_j\text{ converge to positive numbers } \bar\beta,\bar \gamma, \text{ resp.}
\eeqs
Then
\beq\label{doub1}
y_j\le (2A)^{G_{j} \bar \alpha} \max\Big\{  y_0^{\gamma_0\ldots\gamma_{j-1}},y_0^{\beta_0\ldots\beta_{j-1}} \Big\}\quad \forall j\ge 1,
\eeq
where $G_j=\max\{1, \gamma_{m}\gamma_{m+1}\ldots\gamma_{n}:1\le m\le n< j\}$.
Consequently,
\beq\label{doub2}
\limsup_{j\to\infty} y_j\le (2A)^{G \bar\alpha} \max\{y_0^{\bar \gamma},y_0^{\bar \beta}\},
\quad\text{where  }G=\limsup_{j\to\infty} G_j.
\eeq
\end{lemma}
%================================================================%

For iterations below, we will require $$\kappa_*^2<\kappa_f,$$ which is equivalent to
\beq\label{newps}
p_1,p_2<\frac{\sqrt{\kappa_f}+1}{2-a}\quad \text{and}\quad p_3,p_4<\sqrt{\kappa_f}.
\eeq

Therefore, we assume \eqref{newps} for the remainder of this section.
We now fix $\tilde \kappa$ such that 
\beqs%\label{tkapdef}
\kappa_*< \tilde \kappa<\sqrt{\kappa_f}.
\eeqs
For simplicity, we take 
\beq \label{tkapdef}
\tilde \kappa =\sqrt {\frac{\kappa_*^2+\kappa_f}2}. 
\eeq

%=========================================================%
\begin{theorem}\label{LinfU} 
%Assume \eqref{newps}.
Let $\alpha_0$ be a positive number such that 
\beq\label{xal}
\alpha_0 > \max\{2, \alpha_*,\eta_1\} \quad\text{and}\quad 
\alpha_0 \ge \max\Big\{ \frac{\eta_2}{\tilde\kappa -\kappa_*},\frac{a-\delta}{\kappa_f-\tilde\kappa^2}\Big\}.
\eeq

Then there are positive  constants $C$, $\tilde  \mu$, $\tilde \nu$, $\omega_1$, $\omega_2$, $\omega_3$ such that
if  $T>0$ and $\sigma\in (0,1)$ then
\beq\label{Li1}
\|u\|_{L^{\infty}(U\times(\sigma T,T))}\le C\Big(1+\frac1{\sigma T}\Big)^{\omega_1}(1+T)^{\omega_2} \mathcal M_0^{\omega_3} \max\Big\{ \|u\|^{\tilde \mu}_{L^{\tilde \kappa \alpha_0}(U\times(0,T))},\|u\|^{\tilde \nu}_{L^{\tilde\kappa \alpha_0}(U\times(0,T))}\Big\},
\eeq
where $\mathcal M_0$ is defined by \eqref{Mbar}.
\end{theorem}
%=========================================================%
\begin{proof}
%Let $\varep_0=\tilde \kappa-\kappa_*>0$. Then 
We set up to iterate inequality \eqref{bfinest2} in Lemma \ref{GLk}.
Note that $\alpha=\alpha_0$ satisfies \eqref{alp-Large}.
% with $\tilde\kappa$ now given in  \eqref{tkapdef}.
%The choice of $\alpha_0$  \eqref{xal} gives
%\beq
%\kappa(\alpha_0)= \kappa_f  -\frac{a-\delta}{\alpha_0}\ge \tilde \kappa^2.
%\eeq
Define for $j\ge 0$, $\beta_j=\tilde \kappa^j\alpha_0$.
%\beq\label{Betadef}
%\beta_j=\tilde \kappa^j\alpha_0.
%\eeq
Since $\tilde \kappa>1$, the sequence $\{\beta_j\}_{j=0}^\infty$ is increasing, and hence, so is the sequence $\{\kappa(\beta_j)\}_{j=0}^\infty$.
We then have from the last fact and the  choice of $\alpha_0$ that
\beq\label{kk}
\kappa(\beta_j)\ge \kappa(\alpha_0)= \kappa_f  -\frac{a-\delta}{\alpha_0}\ge \tilde \kappa^2.
\eeq

For $j\ge 0$, let $t_j=\sigma T(1-\frac 1{2^j})$.
For $j\ge 0$, applying \eqref{bfinest2} of Theorem \ref{GLk} with $\alpha=\beta_{j}$, $T_2=t_{j+1}$ and $T_1=t_j$, we have 
\begin{align}
\| u\|_{L^{\kappa(\beta_{j})\beta_{j}}(U\times (t_{j+1},T))}
&\le  A_{\beta_{j}}^\frac 1{\beta_{j}}\Big( \| u\|_{L^{\tilde\kappa \beta_{j}}(U\times(t_j,T))}^{\nu_5(\beta_{j})}+\| u\|_{L^{\tilde\kappa \beta_{j}}(U\times(t_j,T))}^{\nu_7(\beta_{j})}\Big)^\frac 1{\beta_{j}} \notag \\
&=  A_{\beta_{j}}^\frac1{\beta_{j}}\Big( \| u\|_{L^{\beta_{j+1}}(U\times(t_j,T))}^{\tilde r_j}+\| u\|_{L^{\beta_{j+1}}(U\times(t_j,T))}^{\tilde s_j}\Big)^\frac 1{\beta_{j}},\label{ukbb}
\end{align}
where
$\tilde r_j=\nu_5(\beta_{j})$ and $\tilde s_j=\nu_7(\beta_{j})$, see formula \eqref{tilrsdef}.
Note from \eqref{kk} that 
$$\kappa(\beta_{j})\beta_{j}\ge \tilde \kappa^2 \beta_{j}=\beta_{j+2}.$$
 
Define for $j\ge 0$ that    $\mathcal Q_j=U\times (t_j,T)$ and
$Y_j=\| u\|_{L^{\beta_{j+1}}(\mathcal Q_j)}.$
By H\"older's inequality
\beq\label{uholder}
 Y_{j+1}= \| u\|_{L^{\beta_{j+2}}(\mathcal Q_{j+1})}\le  |\mathcal Q_{j+1}|^{\frac{1}{\beta_{j+2}}-\frac{1}{\kappa(\beta_{j})\beta_{j}}}  \| u\|_{L^{\kappa(\beta_{j})\beta_{j}}(\mathcal Q_{j+1})}
\le (1+|\mathcal Q_0|)^\frac 1{\beta_{j}}  \| u\|_{L^{\kappa(\beta_{j})\beta_{j}}(\mathcal Q_{j+1})}.
\eeq

Combining inequalities \eqref{ukbb} and \eqref{uholder}  gives
\beq\label{YwithQ}
Y_{j+1}\le \widehat A_j^\frac{1}{\beta_{j}}\big( Y_j^{\tilde{r}_j}+Y_j^{\tilde{s}_j}\big)^{\frac 1{\beta_{j}}} \quad\forall j\ge 0,\quad \text{where }
\widehat A_j=  (1+|\mathcal Q_0|) A_{\beta_{j}}.
\eeq

Now we estimate $\widehat{A}_j$.  
Note that $|\mathcal Q_0|=|Q_T|\le (1+|U|)(1+T)$.
From \eqref{Atildef2} and \eqref{YwithQ},  we have 
\begin{align*}
\widehat A_j
&\le c_7  (1+|U|) \beta_{j}^{6-a}\Big(1+\frac{2^{j+1}}{\sigma T}\Big)^2(1+T)^3 \mathcal M_0^2\\
&\le c_7(1+|U|)(\tilde \kappa^j\alpha_0)^{6-a}4^{j+1}\Big(1+\frac{1}{\sigma T}\Big)^2(1+T)^3 \mathcal M_0^2
\le A_{T,\sigma,\alpha_0}^{j+1},
\end{align*}
where 
\beq\label{AAA}
 A_{T,\sigma,\alpha_0}=\max\big\{ 4\tilde \kappa^{6-a},4c_7(1+|U|) \alpha_0^{6-a}(1+\frac{1}{\sigma T})^2(1+T)^3 \mathcal M_0^2\big\}\ge 1.
 \eeq
Hence 
\beq\label{YY}
Y_{j+1}\le A_{T,\sigma,\alpha_0}^{\frac{j+1}{\beta_j}}\big(Y_j^{\tilde{r}_j}+Y_j^{\tilde{s}_j}\big)^{\frac 1{\beta_j}}.
\eeq

We have
\beqs
\sum_{j=0}^\infty \frac {j+1} { \beta_j}= \frac 1{ \alpha_0}\sum_{j=0}^\infty \frac {j+1} {\tilde \kappa^j}<\infty.
\eeqs

Note that
\beq\label{q12}
1\ge \frac{\tilde r_j}{\beta_j}=1 - \frac{h_1}{\tilde \kappa^j \alpha_0} \quad\text{and}\quad
1\le  \frac{\tilde s_j}{\beta_j}%=\frac{\nu_6}{\beta_j+\delta-a}
= \frac{\beta_j+h_2}{\beta_j+\delta-a}
= 1+\frac{h_2+a-\delta}{\tilde \kappa^j \alpha_0+\delta -a}.
\eeq
Then it is elementary to show that the products 
\beq \label{mnutil}
\tilde \mu =  \prod_{j=0}^\infty \frac{\tilde r_j}{\beta_j}\quad\text{and}\quad \tilde \nu = \prod_{j=0}^\infty \frac{\tilde s_j}{\beta_j}
\eeq 
are  positive numbers. By \eqref{YY} and Lemma \ref{Genn}, we obtain
\beq\label{limY}
\limsup_{j\to\infty} Y_{j}\le (2A_{T,\sigma,\alpha_0})^\omega\max\{Y_0^{\tilde\mu}, Y_0^{\tilde\nu}\},
\eeq
where $\omega=\mathcal G\sum_{j=0}^{\infty}\frac {j+1}{\beta_j}$ with  $\mathcal G=\prod_{k=1}^\infty (\tilde s_k/\beta_k)\in(0,\infty)$.

Note that $Y_0=\| u\|_{L^{\beta_1}(\mathcal Q_0)}$ and, by \eqref{AAA},
\beqs
(2A_{T,\sigma,\varphi})^\omega\le C\Big(1+\frac1{\sigma T}\Big)^{\omega_1}(1+T)^{\omega_2} \mathcal M_0^{\omega_3},
\eeqs
where $\omega_1=2\omega$, $\omega_2=3\omega$ and $\omega_3=2\omega$.
Then estimate \eqref{Li1} follows \eqref{limY}.
\end{proof}
%==================================================================%

\begin{remark}
Estimate \eqref{Li1}  can be rewritten as
$$
\|u\|_{L^{\infty}(U\times(\sigma T,T))}\le C_{\sigma, T,f_1,f_2,\varphi_1,\varphi_2}
\Big( \|u\|^{\tilde \mu}_{L^{\tilde \kappa \alpha_0}(U\times(0,T))}+\|u\|^{\tilde \nu}_{L^{\tilde\kappa \alpha_0}(U\times(0,T))}\Big).
$$
The bound on the right-hand side is quasi-homogeneous in $\|u\|_{L^{\tilde\kappa \alpha_0}(U\times(0,T))}$, hence, relevant for both small and large values of  $\|u\|_{L^{\tilde\kappa \alpha_0}(U\times(0,T))}$.
This is different from the commonly obtained estimates when a positive constant is added to the inequality's right-hand side.
This global (in space) estimate extends our own version in \cite{CHK1}, and the local one in \cite{Surnachev2012}.
\end{remark}

Combining Theorems \ref{lathm} and \ref{LinfU} gives the following particular estimates of in terms of initial and boundary data.
%=============================%
\begin{theorem}\label{thm45}
Let $\alpha_0$ satisfy
\beq\label{alzero}
\alpha_0 > \max\Big\{2, \alpha_*,\eta_1,\frac{\eta_0}{\tilde\kappa}\Big\} \quad\text{and}\quad 
\alpha_0 \ge \max\Big\{ \frac{na}{(1-a)\tilde\kappa},\frac{\eta_2}{\tilde\kappa -\kappa_*},\frac{a-\delta}{\kappa_f-\tilde\kappa^2}\Big\}.
\eeq 

Let $\tilde  \mu$, $\tilde \nu$, $\omega_1$, $\omega_2$, $\omega_3$ be the same constants as in Theorem \ref{LinfU}, and denote $\beta_1=\tilde \kappa\alpha_0$.

\begin{enumerate}
\item[\rm(i)] If $T>0$ satisfies \eqref{Tcond1} for $\alpha=\beta_1$, then for $0<\varepsilon<\min\{1,T\}$, one has 
%\alert{Below, replace $\alpha_0$ with $\beta_1=\tilde \kappa\alpha_0$.}
\beq\label{Lb1}
\|u\|_{L^{\infty}(U\times(\varepsilon,T))}\le C\varepsilon^{-\omega_1}(1+T)^{\omega_2} \mathcal M_0^{\omega_3} \cdot\max\Big\{ \Big(\int_0^T\mathcal{V}(t)dt\Big)^{\frac{\tilde \mu}{\beta_1}}, \Big(\int_0^T\mathcal{V}(t)dt\Big)^{\frac{\tilde \nu}{\beta_1}}\Big\}, 
\eeq
where  $\mathcal M_0$ is defined by \eqref{Mbar}, and for  $0\le t\le T$
\beq
\mathcal{V}(t)=\Big\{\Big(1+\int_Uu_0(x)^{\beta_1} dx\Big)^{-\nu_4(\beta_1)}- \frac 1{C_1(\beta_1)} \int_0^t \Upsilon(\tau)d\tau \Big\}^{-\frac{1}{\nu_4(\beta_1)}}.
\eeq
\item[\rm(ii)] If $T>0$ satisfies \eqref{Tcond2} for $\alpha=\beta_1$, then  for $0<\varepsilon<\min\{1,T\}$, one has
\beq\label{Lb3}
\|u\|_{L^{\infty}(U\times(\varepsilon,T))}\le C\varepsilon^{-\omega_1}(1+T)^{{\omega_2}+\frac{\tilde\nu}{\beta_1}}(1+ \|u_0(x)\|_{L^{\beta_1}(U)})^{\tilde\nu} \mathcal M_0^{\omega_3}.
\eeq
%where $C>0$ depends on $\beta_1$.
\end{enumerate}
%\beq
%\|u(t)\|_{L^\infty}\le bound on initial and boundary data.
%\eeq

Above $C$ is a positive constant independent of $T$ and $\varep$.
\end{theorem}
%==========================%
\begin{proof}
Note that $\alpha_0$ satisfies \eqref{xal}, and $\alpha=\beta_1$ satisfies  \eqref{al3}. The constant $C>0$ in this proof is generic.

(i) Applying estimate \eqref{Li1} to $\sigma T=\varepsilon$ we have 
\beq\label{Li2}
\begin{split}
\|u\|_{L^{\infty}(U\times(\varepsilon,T))}
&\le C\varepsilon^{-\omega_1}(1+T)^{\omega_2} \mathcal M_0^{\omega_3}\\
&\quad\cdot  \max\Big\{ \Big(\int_0^T\int_U|u|^{\beta_1}dxdt\Big)^{\frac{\tilde \mu}{\beta_1}},
 \Big(\int_0^T\int_U|u|^{\beta_1}dxdt\Big)^{\frac{\tilde \nu}{\beta_1}}\Big\}.
\end{split}
\eeq
By \eqref{La}, we have 
\beq\label{Lb}
\int_U u^{\beta_1}(x,t)dx \le  \mathcal{V}(t).
\eeq 
Then \eqref{Li2} and \eqref{Lb} yield \eqref{Lb1}.

(ii) If $T$ satisfies \eqref{Tcond2} with $\alpha=\beta_1$, then by \eqref{La2}, we have 
\beq\label{Lb2}
\int_U u^{\beta_1}(x,t)dx \le 2\Big(1+\int_U u_0^{\beta_1}(x)dx\Big)\quad\text{for all}\quad t\in[0,T].
\eeq 
Combining \eqref{Li2} with \eqref{Lb2}, and noticing that  $ \tilde\mu < \tilde\nu$, we obtain 
\beqs%\label{Li3}
\begin{split}
&\|u\|_{L^{\infty}(U\times(\varepsilon,T))}
\le C\varepsilon^{-\omega_1}(1+T)^{\omega_2} \mathcal M_0^{\omega_3} \\
&\quad \cdot \max\Big\{ \Big(2\int_0^T\Big(1+\int_U u_0^{\beta_1}(x)dx\Big)dt\Big)^{\frac{\tilde \mu}{\beta_1}},
 \Big(2\int_0^T\Big(1+\int_U u_0^{\beta_1}(x)dx\Big)dt\Big)^{\frac{\tilde \nu}{\beta_1}}\Big\}\\
 &\le C\varepsilon^{-\omega_1}(1+T)^{\omega_2} \mathcal M_0^{\omega_3} \cdot (1+T)^{\frac{\tilde \nu}{\beta_1}}
\Big(1+\int_U u_0^{\beta_1}(x)dx\Big)^{\frac{\tilde \nu}{\beta_1}}, 
 \end{split}
\eeqs
and estimate \eqref{Lb3} follows. 
\end{proof}

%==================================================================%

%%%%%%%%%%%%%%%%%%%%%%%%%%%%%%%%%%%%%%%%%%%%%%%%%%%%%%%%%%%%%%%%%%%%%
\myclearpage
\section{Gradient estimates}\label{gradsec}
%%%%%%%%%%%%%%%%%%%%%%%%%%%%%%%%%%%%%%%%%%%%%%%%%%%%%%%%%%%%%%%%%%%%%
In Theorem \ref{lathm} we have an estimate for $\int_0^T \int_U |\nabla u(x,t)|^{2-a}u(x,t)^{\alpha-\lambda-1} dxdt$. 
Note that the double integral is in both spatial and time variables, and it is not yet a direct estimate for the gradient.
Therefore, we focus, in the following,  on estimating $\int_U |\nabla u(x,t)|^{2-a}dx$ directly for $t>0$. 

In connection with properties  \eqref{Zprop} and \eqref{Bprop} in Assumption (A1), it is natural to make the following assumptions. 
%These assumptions are natural growth conditions considering \eqref{Zprop} and \eqref{Bprop}.

\bigskip
\noindent\textbf{Assumption (A2).} 
The function $Z(u)$ satisfies 
\beq\label{Zp}
|Z'(u)|\le d_4 u^{\ellz-1} \quad\forall u\in [0,\infty),
\eeq
for some constant $d_4>0$, and there are non-negative functions $\varphi_3(x,t)$ and $\varphi_4(x,t)$  defined on $\Gamma\times(0,\infty)$ such that
\beq\label{Bt}
\Big|\frac{\partial B(x,t,u)}{\partial t}\Big|\le \varphi_3(x,t) +\varphi_4(x,t) u^\ellb.
\eeq
We also assume a slightly stronger version of \eqref{Bprop} and \eqref{fprop}, namely,
\beq \label{Babs}
 |B(x,t,u)|\le \varphi_1(x,t)+\varphi_2(x,t)u^\ellb,
\eeq
\beq\label{fabs}
|f(x,t,u)|\le f_1(x,t)+f_2(x,t)u^\ellf.
\eeq

\bigskip
To deal with the boundary condition, we define for any $x\in\Gamma$, $t\ge 0$, $u\ge 0$
\beq\label{Qdef}
Q(x,t,u)=\int_0^u B(x,t,v)dv.
\eeq

Then by Assumptions (A1) and (A2), we have
\beq\label{Q1}
|Q(x,t,u)|\le \varphi_1(x,t) u+\varphi_2(x,t) u^{\ellb+1}
\eeq
and 
\beq\label{Q2}
\Big|\frac{\partial Q(x,t,u)} {\partial t}\Big|\le \varphi_3(x,t) u +C\varphi_4(x,t) u^{\ellb+1} .
\eeq

\bigskip
\noindent\textbf{Assumption (A3).}
 \beq\label{lamel}
2\ellz >\lambda+1.
\eeq

For our original problem \eqref{maineq}, $\ellz=2\lambda$, then condition \eqref{lamel} becomes $\lambda>1/3$.
For slightly compressible fluids, $\lambda=1$.
For ideal gases, $\lambda=1/2$.
For isentropic gas flows, from the data in \cite{GasBook2001} or section III of \cite{LangePierrus2000}, all values of the specific heat ratio $\gamma$ belong to the interval $(1,2)$, therefore $\lambda=1/(1+\gamma)$, see \eqref{gaspara}, satisfies $1/3<\lambda< 1/2$. Thus \eqref{lamel} is naturally met in all cases.

\bigskip
For the gradient estimates, same as in \cite{ABHI1,HI1}, we will make use of the following function
\beqs
H(\xi)=\int_0^{\xi^2} K(\sqrt{s}) ds \quad \text{for } \xi\ge 0. 
\eeqs
The function $H(\xi)$ satisfies
\beq
K(\xi)\xi^2 \le H(\xi)\le 2K(\xi)\xi^2, \label{i:ineq4}
\eeq
and hence, as a consequence of \eqref{i:ineq2} and \eqref{i:ineq4}, 
\beq
d_3(\xi^{2-a}-1) \le H(\xi) \le 2 d_2\xi^{2-a}. \label{i:ineq5}
\eeq

Based on the structure of the PDE in \eqref{upb}, we define an intermediate quantity
\beqs
\mathcal I(t)=\int_U H\Big(|\nabla u(x,t)+Z(u(x,t))|\Big) dx\quad\text{for } t\ge 0,
\eeqs
and for initial values
\beqs
\mathcal Z_0=\int_U u_0^{\lambda+1}(x) dx+\mathcal I(0)\\
+ \int_\Gamma  \Big(\varphi_1(x,0)u_0(x)+ \varphi_2(x,0)u_0^{\ellb+1}(x) \Big) d\sigma.
\eeqs

First, we estimate $\mathcal I(t)$ in terms of certain Lebesgue norms of $u(x,t)$.

\begin{proposition}\label{prop51}
 For $t> 0$,
\begin{multline}\label{H1}
\mathcal I(t)\le   2\mathcal Z_0 +  C\Big\{ t+1  + \int_U u^{\eta_3}(x,t)dx+   \int_0^t \int_U u^{\eta_4} (x,\tau)dxd\tau
+ N_1(t)+ \int_0^t N_2(\tau) d\tau\Big\},
\end{multline}
where $C$ is a positive constant,
\begin{align*}
\eta_3&=\max\Big\{\frac{(2-a)(2\ellb+1)}{1-a}, (2-a)\ellz \Big\},\\
\eta_4&= \max\Big\{\frac{(2-a)(2\ellb+1)}{1-a},  \frac{2\ellz(2-a)}{a},2(2\ellf+1) \Big\},
\end{align*}
and
\begin{align*}
 N_1(t)&=\int_\Gamma \Big(\varphi_1^{\eta_5}(x,t)+\varphi_2^2(x,t)\Big)d\sigma,\\
N_2(t) &= N_1(t)+\int_\Gamma \Big(\varphi_3^{\eta_5}(x,t)+\varphi_4^2(x,t)\Big)d\sigma + \int_U \Big(f_1^{\eta_6}(x,t)+f_2^4(x,t)\Big)dx
 \end{align*}
 with 
 \beqs
 \eta_5=\frac{2(1+\ellb)}{1+2\ellb},\quad
 \eta_6=\max\Big\{ \frac{2(1+\ellf)}{1+2\ellf} ,\frac{4(2\ellf+1-\lambda)}{4\ellf+1-\lambda}\Big\}.
 \eeqs
\end{proposition}
\begin{proof}
The proof is divided into three steps below. As usual, constants $C$, $C_\varep$ are positive and generic.

\textbf{Step 1.} Multiplying the PDE in \eqref{upb} by $u$, integrating over $U$ and using integration by parts, we have
 \begin{align*}
\frac\lambda{\lambda+1}\ddt \int_U u^{\lambda+1}dx&= \int_U \nabla\cdot \left(K(|\nabla u+Z(u)|)(\nabla u+Z(u)) \right) u dx +\int_U f u dx\\
&= - \int_U K(|\nabla u+Z(u)|)(\nabla u+Z(u))\cdot \nabla u  dx 
+\int_{\Gamma} B u  d\sigma +\int_U f u dx\\
&= - \int_U K(|\nabla u+Z(u)|)|\nabla u+Z(u)|^2 dx \\
&\quad +\int_U K(|\nabla u+Z(u)|)(\nabla u+Z(u))\cdot Z(u)  +\int_\Gamma  B u d\sigma +\int_U f u dx.
\end{align*}

On the right-hand side of the last equation, using relation \eqref{i:ineq4} for the first integral, relations \eqref{i:ineq1}, \eqref{Zprop} for the second integral, \eqref{Bprop} for the third, and \eqref{fprop} for the last, we have 
\beq\label{predH}
\begin{aligned}
\frac\lambda{\lambda+1}\ddt \int_U u^{\lambda+1}dx&\le - \frac 12 \int_U H(|\nabla u+Z(u)|)dx  +C \int_U |\nabla u+Z(u)|^{1-a} u^{\ellz}dx\\
&\quad +\int_\Gamma (\varphi_1+\varphi_2 u^\ellb) u d\sigma+\int_U (f_1+f_2 u^\ellf) udx\\
& \eqdef -\frac12 \mathcal I + J_1+J_2+J_3.
\end{aligned}
\eeq

In the following,  $\varep\in(0,1)$.

$\bullet$ For $J_1$, using Young's inequality with powers $\frac{2-a}{1-a}$ and $2-a$, and then relation \eqref{i:ineq5},
\beq\label{J2}
J_1\le \varep \mathcal I + C_\varep+C_\varep\int_U u^{\ellz(2-a)} dx.
\eeq
%Then 
%\beq\label{Jmid}
%\frac\lambda{\lambda+1}\ddt \int_U u^{\lambda+1}dx\le - \frac 14 \mathcal I +C+C\int_U u^{\ellz(2-a)} dx +J_4+J_5.
%\eeq

$\bullet$ For $J_2$, using Young's inequality, we have 
\begin{align*}
J_2\le 2 \|u^{1+\ellb}\|^2_{L^{2}(\Gamma)}+ N_1(t).
\end{align*}
Next, applying inequality \eqref{etrace10} for %$\alpha=s=2$, $p=2-a$, and  
$\alpha =s=2(\ellb+1)$, $p=2-a$ to  the boundary integral $\|u^{1+\ellb}\|^2_{L^{2}(\Gamma)}$, we have 
\begin{align*}
 J_2
&\le   \varep \int_U |\nabla u|^{2-a} dx + C\|u\|_{L^{2(\ellb+1)}}^{2(\ellb+1)}+C_\varep\|u\|_{L^{s_1}}^{s_1}+ N_1(t)
\end{align*}
with %$\mu_8 =\frac{2-a}{1-a}$ and 
$$s_1 =2(\ellb+1)+\frac{2\ellb+a}{1-a}=\frac{(2-a)(2\ellb+1)}{1-a}.$$
By \eqref{ee3} and then relations \eqref{Zprop}, \eqref{i:ineq5}, 
\begin{align*}
  \int_U |\nabla u|^{2-a} dx &\le 2\int_U |\nabla u + Z(u)|^{2-a} dx +2\int_U |Z(u)|^{2-a} dx\\
  &\le \frac{2}{d_3}\mathcal I +C +C\int_U u^{(2-a)\ellz} dx.   
\end{align*}
Thus
\beqs
 J_2\le  \frac{2\varep}{d_3} \mathcal I  +C+C\int_U u^{(2-a)\ellz} dx + C\int_U u^{2(\ellb+1)}+C_\varep\int_U u^{s_1} dx+N_1(t).
\eeqs
Since $(2-a)\ellz, 2(\ellb+1),s_1\le \eta_3$, by using Young's inequality, we have 
\beq\label{J4}
\begin{aligned}
 J_2
&\le   \frac{2\varep}{d_3}  \mathcal I  +C_\varepsilon+C_\varepsilon\int_U u^{\eta_3} dx + N_1(t).
 \end{aligned}
\eeq
%where $$\eta_3=\max\{2(\ellb+1),s_1, 2(\ellz+1) \}.$$
$\bullet$ For $J_3$, using Young's inequality
\beq\label{J5}
J_3\le 2\int_U u^{2(\ellf+1)}  dx+  \int_U \Big(f_1^\frac{2(1+\ellf)}{1+2\ellf}+f_2^4\Big)dx.
\eeq
Combining \eqref{predH}, \eqref{J2}, \eqref{J4} and \eqref{J5} gives %using Young inequality give
\beqs
\begin{split}
&\frac\lambda{\lambda+1}\ddt \int_U u^{\lambda+1}dx + \Big(\frac 1 2 -\varep - \frac{2\varep}{d_3} \Big ) \mathcal I \\
&\le C_\varep\Big( 1+\int_U (u^{\ellz(2-a)}+u^{\eta_3} +u^{2(\ellf+1)} ) dx\Big)+N_3(t),
\end{split}
\eeqs
where
\beqs
N_3(t)= N_1(t)+\int_U \Big(f_1^\frac{2(1+\ellf)}{1+2\ellf}(x,t)+f_2^4(x,t)\Big)dx.
\eeqs
Taking $\varep$ sufficiently small and using Young's inequality show that
\beq\label{dif1}
\ddt \int_U u^{\lambda+1}dx  + \mathcal I 
\le C + C \int_U u^{s_2}  dx + CN_3(t).
\eeq
where
\beqs
s_2 =\max\{  2(\ellf+1), \eta_3  \}.
\eeqs 

\textbf{Step 2.} Multiplying both sides of the PDE in \eqref{upb} by  $ u_t $, integrating over $U$ and using the boundary condition, one has   
\[
\lambda \int_U u^{\lambda-1} u_t^2 dx +\int_U  K (|\nabla u+Z(u)|)(\nabla u+Z(u))\cdot \nabla u_t dx  = \int_\Gamma  B  u_t  d\sigma+\int_U f u_t dx.
\]
Then by \eqref{i:ineq1}, we have
\begin{align*}
&\lambda \int_U u^{\lambda-1} u_t^2 dx +\int_U  K (|\nabla u+Z(u)|)(\nabla u+Z(u))\cdot (\nabla u+Z(u))_t dx  \\
&= \int_U  K (|\nabla u+Z(u)|)(\nabla u+Z(u))\cdot (Z(u))_t dx  + \int_\Gamma  B  u_t  d\sigma + \int_U f u_t dx\\
&\le  C \int_U  |\nabla u+Z(u)|^{1-a} |Z'(u)||u_t| dx + \int_\Gamma  B  u_t  d\sigma + \int_U f u_t dx.
\end{align*}

Hence by definition of $H(\xi)$, we have 
\beq\label{dif5} 
\begin{split}
&\lambda \int_U u^{\lambda-1} u_t^2 dx +\frac12\ddt \int_U  H (|\nabla u+Z(u)|) dx  \\
&\le  C \int_U  |\nabla u+Z(u)|^{1-a} |Z'(u)||u_t| dx + \int_\Gamma  B  u_t  d\sigma + \int_U f u_t dx.
\end{split}
\eeq

$\bullet $ We bound the first term of the right hand side in \eqref{dif5} by  \eqref{Zp} and then Young's inequality,
\begin{align*}
& C \int_U  |\nabla u+Z(u)|^{1-a}  |Z'(u)||u_t| dx 
\le C \int_U  |\nabla u+Z(u)|^{1-a} u^{\ellz-1} |u_t| dx\\
&\le \varep \int_U  |\nabla u+Z(u)|^{2-a}  dx + C_\varep \int_U  u^{(2-a)(\ellz-1)} |u_t|^{2-a} dx.
\end{align*}
By \eqref{i:ineq5} and then using Young's inequality one more time  on the last term, we have 
\beq\label{Z1}
C \int_U  |\nabla u+Z(u)|^{1-a}  |Z'(u)||u_t| dx \le \frac{\varep}{d_3} \mathcal I  + C+\frac\lambda4 \int_U u^{\lambda-1} u_t^2 dx+C_\varep \int_U  u^\frac{(2-a)(2\ellz-\lambda-1)}{a}  dx.
\eeq

$\bullet $  For the middle term on the right-hand side of \eqref{dif5}, using definition \eqref{Qdef} we have
\beqs
B(x,t,u(x,t))u_t(x,t) = \ddt Q(x,t,u(x,t)) - \frac{\partial Q} {\partial t}\Big|_{(x,t,u(x,t))}.
\eeqs
Thus \eqref{Bt} yields
\beq\label{Btreln}
\begin{aligned}
 \int_\Gamma  B  u_t  d\sigma 
&=\ddt \int_\Gamma Q(x,t,u)d\sigma - \int_\Gamma \frac{\partial Q(x,t,u)} {\partial t} d\sigma\\
&\le  \ddt \int_\Gamma Q(x,t,u)d\sigma +C\int_\Gamma  (\varphi_3 u +\varphi_4 u^{\ellb+1})d\sigma.
\end{aligned}
\eeq
For the last term of \eqref{Btreln}, similar to estimate \eqref{J4} we have  
\beq\label{J40}
\begin{aligned}
 C\int_\Gamma  (\varphi_3 u +\varphi_4 u^{\ellb+1})d\sigma
&\le  \frac{2\varep}{d_3} \mathcal I  + C_\varepsilon+C_\varepsilon\int_U u^{\eta_3} dx 
+C\int_\Gamma\Big( \varphi_3^{\eta_5}+\varphi_4^2\Big)d\sigma.
 \end{aligned}
\eeq

$\bullet $  For the last term of \eqref{dif5}, using Cauchy's inequality, \eqref{fabs} and Young's inequality, we  obtain   
\beq\label{fut}
\begin{aligned}
 \int_U f u_t dx &\le \frac\lambda4\int_U u^{\lambda-1} u_t^2 dx + C \int_U f^2u^{1-\lambda} dx\\
&\le \frac\lambda4 \int_U u^{\lambda-1} u_t^2 dx + C \int_U ( f_1^2+f_2^2u^{2\ellf})u^{1-\lambda} dx\\
&\le \frac\lambda4 \int_U u^{\lambda-1} u_t^2 dx +C\int_Uu^{2(2\ellf+1-\lambda)} dx 
+C \int_U\Big( f_1^{\frac{4(2\ellf+1-\lambda)}{4\ellf+1-\lambda}}+f_2^4\Big)dx.
 %C (\norm{f_1^2}^2+ \norm{f_2^2}^2 ) + \int_U u^{2(\lambda-1)}+u^{2(2\ellf+1-\lambda)} dx   .
\end{aligned}
\eeq

 Then combining  \eqref{dif5},  \eqref{Z1},  \eqref{Btreln}, \eqref{J40}, and \eqref{fut},  we obtain
 \beq\label{dif6} \begin{split}
&\frac\lambda2\int_U u^{\lambda-1} u_t^2 dx +\frac12\ddt \mathcal I   - \ddt \int_\Gamma Q(x,t,u)d\sigma \\
&\le  \frac{3\varep}{d_3}\mathcal I  + C_\varepsilon +C_\varepsilon\int_U \Big( u^\frac{(2-a)(2\ellz-\lambda-1)}{a}  +u^{\eta_3}+ u^{2(2\ellf+1-\lambda)}\Big )dx  + CN_4(t),
\end{split}
\eeq
where
 \beq
 N_4(t)=\int_\Gamma\Big( \varphi_3^{\eta_5}+\varphi_4^2\Big)d\sigma +\int_U\Big( f_1^{\frac{4(2\ellf+1-\lambda)}{4\ellf+1-\lambda}}+f_2^4\Big)dx.
\eeq

%Let $C_{u,T}$ denote the bound depending on $u$, $f_i$, $\varphi_i$ and $T$.
\textbf{Step 3.} Let $t\in (0,T)$. Note by assumption \eqref{lamel} that the power $(2-a)(2\ellz-\lambda-1)/a$ in \eqref{dif6} is positive.
%Note that
%\beq
%\eta_3,s_2,\frac{(2-a)(2\ellz-\lambda-1)}{a},2(2\ellf+1-\lambda)\le s_3.
%\eeq
Summing \eqref{dif1}, \eqref{dif6}  with sufficiently small $\varep$, and using Young's inequalities, we have
\beq\label{dNs}
\frac d {dt}\Big(\int_U u^{\lambda+1} dx+ \mathcal I  - \int_\Gamma  Q(x,t,u)  d\sigma\Big)% + C\mathcal I \\
 \le C+ C \int_U u^{s_3} dx  +C (N_3(t)+N_4(t))
 \eeq 
with 
\beqs
s_3 =\max\Big\{ s_2, \frac{(2-a)(2\ellz-\lambda-1)}{a}, 2(2\ellf+1-\lambda) \Big\}.
\eeqs
%\begin{align*}
%N_3(t) &= \int_\Gamma \Big(\varphi_1(x,t)^{\frac{2(1+\ellb)}{1+2\ellb}}+\varphi_2(x,t)^2+\varphi_3(x,t)^{\frac{2(1+\ellb)}{1+2\ellb}}+\varphi_4(x,t)^2\Big)d\sigma\\
%&\quad + \int_U \Big(f_1(x,t)^\frac{2(1+\ellf)}{1+2\ellf}  +f_1(x,t)^{\frac{4(2\ellf+1-\lambda)}{4\ellf+1-\lambda}}+f_2(x,t)^4\Big)dx.
% \end{align*}
 
For simplicity, we bound $s_3\le \eta_4$, and  $u^{s_3}\le 1+u^{\eta_4}$. 
Also, the powers $\frac{2(1+\ellf)}{1+2\ellf}$ and $\frac{4(2\ellf+1-\lambda)}{4\ellf+1-\lambda}$ of $f_1$ in $N_3(t)$ and $N_4(t)$ are less than or equal to $\eta_6$, then by Young's inequality
$$N_3(t)+N_4(t)\le C(1+N_2(t)).$$
Therefore, \eqref{dNs} yields
\beq\label{H:ineq}
\frac d {dt}\Big(\int_U u^{\lambda+1} dx+ \mathcal I  - \int_\Gamma  Q(x,t,u)  d\sigma\Big)
 \le C+ C \int_U u^{\eta_4} dx  +C N_2(t).
\eeq 

Integrating \eqref{H:ineq} in time from $0$ to $t$ gives
\begin{multline}\label{ineqb}
\int_U u^{\lambda+1}(x,t) dx+\mathcal I (t)  
\le \int_U u_0(x)^{\lambda+1} dx+ \mathcal I (0) - \int_\Gamma  Q(x,0,u_0(x))  d\sigma + \int_\Gamma  Q(x,t,u)  d\sigma\\
+ Ct+ C \int_0^t \int_U u^{\eta_4} dxd\tau  +C \int_0^t N_2(\tau) d\tau .
\end{multline}
We will neglect the first term on the left-hand side of \eqref{ineqb}. On the right-hand side, by \eqref{Q1} we have
$$ \int_\Gamma  |Q(x,0,u_0(x))|  d\sigma \le \int_\Gamma (\varphi_1(x,0) u_0(x)+\varphi_2(x,0) u_0^{\ellb+1}(x))d\sigma,$$
while for $t>0$, using also estimate \eqref{J4} we obtain %and Young's inequality we obtain
\beq\label{ineqq1}
 \int_\Gamma  Q(x,t,u)  d\sigma \le \int_\Gamma (\varphi_1 u+\varphi_2 u^{\ellb+1})d\sigma\\
 \le  \frac{2\varep}{d_3} \mathcal I  + C_\varepsilon + C_\varep\int_U u^{\eta_3}dx+ CN_1(t). 
\eeq
Again, choosing $\varep$ sufficiently small, we derive from \eqref{ineqb} and \eqref{ineqq1} that      
\beqs%\label{H:est}
%\begin{split}
\frac12\mathcal I (t) \le \mathcal Z_0+ C(t+1) + C \int_0^t \int_U u^{\eta_4} dxd\tau  +C \int_0^t N_2(\tau) d\tau + C\int_U u^{\eta_3}dx+ C N_1(t),
%\end{split}
\eeqs 
and, hence, estimate  \eqref{H1} follows.
\end{proof}

Finally, we combine the estimate in Proposition \ref{prop51} with those in Theorem \ref{lathm}.

\begin{theorem}\label{thm52}
Let
\beqs
\eta_7=\max\Big \{\frac{na}{1-a},\eta_0+1,\eta_4\Big\} \quad\text{and}\quad \eta_8=\nu_4(\eta_7).
\eeqs
\begin{enumerate}
\item[\rm (i)] If $T>0$ satisfies \eqref{Tcond1} for $\alpha=\eta_7$, then for all $t\in(0,T]$
\beq\label{G:est}
\int_U |\nabla u(x,t)|^{2-a} dx  \le C\Big( \mathcal Z_0 +(t+1)(1+\mathcal V_*(t))+ N_1(t)+ \int_0^t N_2(\tau) d\tau\Big ),
\eeq
where
\beqs
\mathcal V_*(t)=  \Big\{\Big(1+\int_Uu_0(x)^{\eta_7} dx\Big)^{-\eta_8}- \frac 1{C_{1}(\eta_7)} \int_0^t (1+\Upsilon(\tau)) d\tau \Big\}^{-\frac{1}{\eta_8}}. 
\eeqs
% where $C>0$ depends on $\alpha$, and
%\beqs
%A(t)=(t+1)(1+\mathcal V_*(t))+ \int_0^t N_2(\tau) d\tau+ N_1(t).
%\eeqs

\item[\rm (ii)] Moreover, if $T>0$ satisfies \eqref{Tcond2} for $\alpha=\eta_7$ then for all $t\in(0,T]$
\beq\label{Gsmall}
\int_U |\nabla u(x,t)|^{2-a} dx  \le C\Big\{ \mathcal Z_0 +(t+1)\Big(1+\int_U u_0^{\eta_7}(x)dx\Big)+ N_1(t)+ \int_0^t N_2(\tau) d\tau\Big\}.
\eeq
\end{enumerate}

Above, $C$ is a positive constant.
\end{theorem}
\begin{proof}
Thanks to \eqref{ee3}, \eqref {i:ineq5} and \eqref{H1}, we have  
 \beqs
\begin{split}
&\int_U |\nabla u|^{2-a} dx \le 2\int_U (|\nabla u+Z(u)|^{2-a}+|Z(u)|^{2-a}) dx 
%\le 2\int_U |\nabla u+Z(u)|^{2-a} dx +C\int_U u^{\ellz(2-a)}dx\\
\le C\mathcal{I}(t) +C+C\int_U u^{\ellz(2-a)}dx\\
&\le C\Big\{ \mathcal Z_0 +  (t+1) + \int_U (u^{\eta_3}+u^{\ellz(2-a)})dx+    \int_0^t \int_U u(x,\tau)^{\eta_4} dxd\tau 
+ N_1(t) +\int_0^t N_2(\tau) d\tau\Big\}.
\end{split}
\eeqs
Using the fact   $\eta_3,\ellz(2-a),\eta_4\le \eta_7$ and \eqref{ee5} to bound $u^{\eta_3},u^{\ellz(2-a)},u^{\eta_4}\le 1+u^{\eta_7}$,
we obtain
\beq\label{wlz1}
\begin{split}
\int_U |\nabla u|^{2-a} dx 
&\le  C\Big\{ \mathcal Z_0 + (t+1)+ \int_U u(x,t)^{\eta_7}dx+     \int_0^t \int_U u(x,\tau)^{\eta_7} dxd\tau  \\
&\quad + N_1(t)+ \int_0^t N_2(\tau) d\tau\Big\}.
\end{split}
\eeq

(i) Note that  $\alpha=\eta_7$ satisfies \eqref{al3}. Let $T>0$ satisfy \eqref{Tcond1} for $\alpha=\eta_7$.
Now using \eqref{La} with $\alpha=\eta_7$, we have
\beq\label{Lai}
\int_U u(x,t)^{\eta_7} dx \le  \mathcal V_*(t).  
\eeq

By \eqref{Lai} and the fact $\mathcal V_*(t)$ is increasing in $t$, we have 
\beq\label{finalineq}
%\begin{split}
\int_U |\nabla u|^{2-a} dx  \le C\Big\{\mathcal Z_0 + (t+1)+   (t+1)\mathcal V_*(t)+N_1(t)+ \int_0^t N_2(\tau) d\tau\Big\},
%\end{split}
\eeq
and, therefore, obtain \eqref{G:est}.   

(ii) Now, assume $T>0$ satisfies \eqref{Tcond1} for $\alpha=\eta_7$. By using \eqref{La2} instead of \eqref{La}, we can replace $\mathcal V_*(t)$ with $2(1+\int_U u_0^{\eta_7}(x)dx)$ in \eqref{Lai} and  \eqref{finalineq}, and consequently, estimate \eqref{Gsmall} follows. \end{proof}

\textbf{Acknowledgement.}  L.~H. acknowledges the support by NSF grant DMS-1412796.

%%%%%%%%%%%%%%%%%%%%%%%%%%%%%%%%%%%%%%%%%%%%%%%%%%%%%%
% Bibliography using BibTeX
%%%%%%%%%%%%%%%%%%%%%%%%%%%%%%%%%%%%%%%%%%%%%%%%%%%%%%
\myclearpage

%\bibliography{paperbaseall}{}

%\bibliographystyle{plain}
\bibliographystyle{acm}
\end{document}